\documentclass[12pt]{amsart}
\usepackage{amsthm}
\usepackage{amsmath}
\usepackage{amssymb} 
\usepackage{graphicx} 
\usepackage{float}
\usepackage{amsfonts}
\usepackage{tikz}
\usepackage{xcolor}
\usepackage{url}
\usepackage{soul}
\usepackage{enumitem}
\usepackage{ulem}

\newtheorem{theorem}{Theorem}[section]
\newtheorem{corollary}[theorem]{Corollary}
\newtheorem{lemma}[theorem]{Lemma}

\newtheorem{remark}[theorem]{Remark}

\DeclareMathOperator{\area}{area}
\DeclareMathOperator{\perimeter}{perim}
\DeclareMathOperator{\height}{height}
\DeclareMathOperator{\width}{width}

\DeclareMathOperator{\recarea}{recarea}

\newcommand{\Z}{\mathbb{Z}}

\begin{document}

\title{Area, Perimeter, Height, and Width of Rectangle Visibility Graphs}

\author[Caughman]{John S. Caughman}\address{Department of Mathematics \& Statistics, Portland State University, Box 751, Portland, OR 97202.}

\author[Dunn]{Charles L. Dunn}\address{Department of Mathematics \& Computer Science, Linfield University, 900 SE Baker Street, McMinnville, OR 97128.}

\author[Laison]{Joshua D. Laison}\address{Department of Mathematics, Willamette University, Salem, OR 97301.}

\author[Neudauer]{Nancy Ann Neudauer}\address{Department of Mathematics \& Computer Science, Pacific University, Forest Grove, OR 97116.}

\author[Starr]{Colin L. Starr}\address{Department of Mathematics, Willamette University, Salem, OR 97301.}

\date{\today}

\begin{abstract}
A rectangle visibility graph (RVG) is represented by assigning to each vertex a rectangle in the plane with horizontal and vertical sides in such a way that edges in the graph correspond to unobstructed horizontal and vertical lines of sight between their corresponding rectangles.  To discretize, we consider only rectangles whose corners have integer coordinates. For any given RVG, we seek a representation with smallest bounding box as measured by its area, perimeter, height, or width (height 
is assumed not to exceed width).

We derive a number of results regarding these parameters. Using these results, we show that these four measures 
are distinct, in the sense that there exist graphs $G_1$ and $G_2$ with $\area(G_1)<\area(G_2)$ but $\perimeter(G_2)<\perimeter(G_1)$, and analogously for all other pairs of these parameters. We further show that there exists a graph $G_3$ with 
representations $S_1$ and $S_2$ such that $\area(G_3)=\area(S_1)<\area(S_2)$ but $\perimeter(G_3)=\perimeter(S_2)<\perimeter(S_1)$. In other words, $G_3$ requires distinct 
representations to minimize area and perimeter. Similarly, such graphs exist to demonstrate the independence of all other pairs of these parameters.

Among graphs with $n \leq 6$ vertices, the empty graph $E_n$ requires largest area. But for graphs with $n=7$ and $n=8$ vertices, we show that the complete graphs $K_7$ and $K_8$  require larger area than $E_7$ and $E_8$, respectively. Using this, we show that for all $n \geq 8$, the empty graph $E_n$ does not have largest area, perimeter, height, or width among all RVGs on $n$ vertices. 
\end{abstract}
\maketitle

\section{Introduction}

\begin{figure}[ht!]
\centering
\includegraphics[width=.6\textwidth]{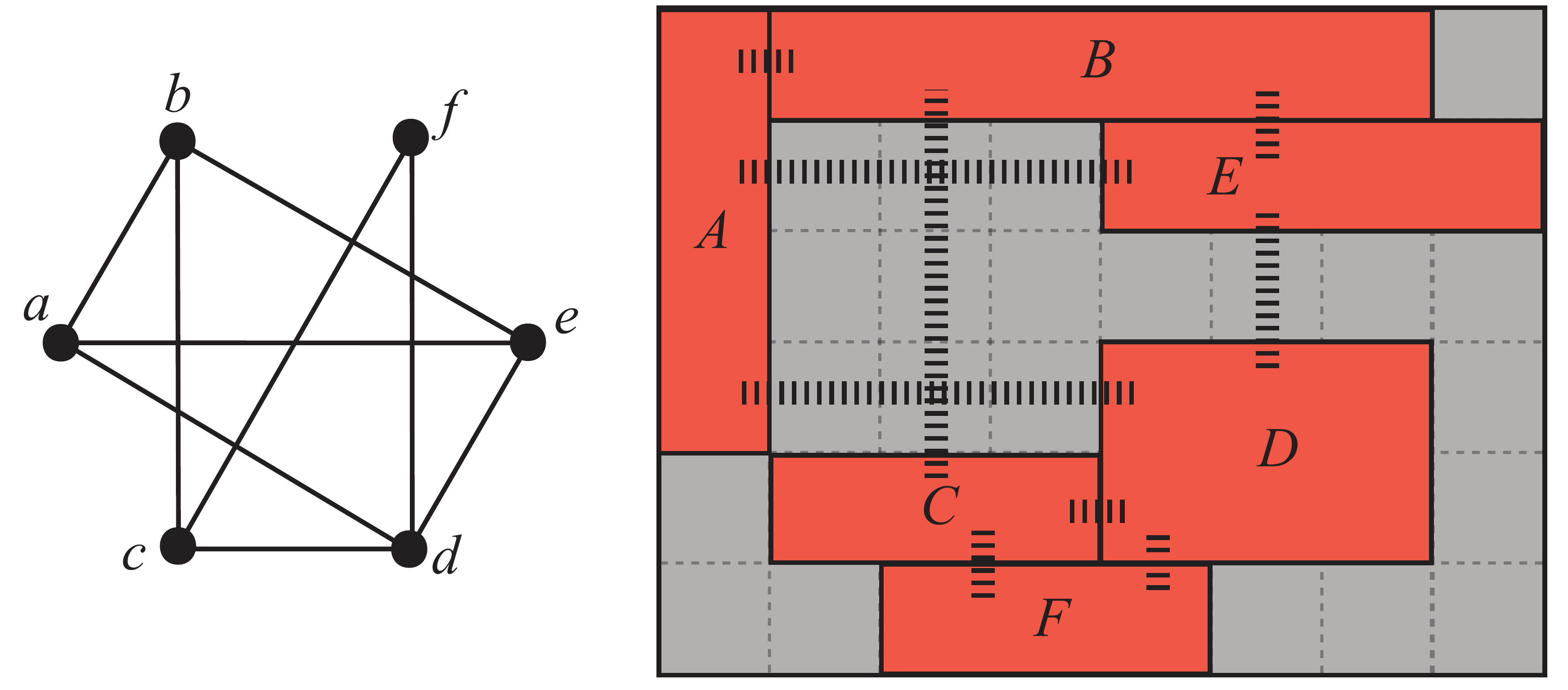}
\caption{A rectangle visibility graph and a corresponding RV-representation with integer rectangles.}
\label{fig:RVG-example} 
\end{figure}

Let $S$ be a set of rectangles in the plane, with vertical and horizontal sides, whose interiors do not intersect.  We say that two rectangles $A$ and $B$ in $S$ \textit{\textbf{see each other}} if there is a vertical or horizontal line segment intersecting the interiors of both $A$ and $B$ and intersecting no other (closed) rectangles in $S$, like the dotted lines in Figure~\ref{fig:RVG-example}. We refer to such segments as \textit{\textbf{lines of sight}}, and under this definition we may consider them to have small positive width. For example, there is no line of sight between rectangles $B$ and $F$ in Figure~\ref{fig:RVG-example}, since a line of sight needs positive width.  

We construct a graph $G$ with a vertex for every rectangle in $S$, and an edge between two vertices if and only if their corresponding rectangles see each other.  We say that $S$ is a \textbf{\textit{rectangle visibility}} (\textit{\textbf{RV-}}) \textit{\textbf{representation}} of $G$, and $G$ is a \textit{\textbf{rectangle visibility graph}} or \textit{\textbf{RVG}}.  We allow rectangles in $S$ to share edges. 

A similar notion of rectangle visibility graph was first introduced in 1976 by Garey, Johnson and So \cite{garey76} as a tool to study the design of printed circuit boards.  Their RVGs have only $1\times 1$ squares, located at a set of lattice points in a grid.  Hutchinson continued work on this problem in 1993 \cite{hutchinson93}.  In \cite{bose1996rectangle}, Bose, Dean, Hutchinson and Shermer described the problem of \textit{two-layer routing} in ``very large-scale integration'' (VLSI) design as follows: 
\begin{quote} {\small
In two-layer routing, one embeds processing components and their connections (sometimes called \textit{wires}) in two layers of silicon (or other VLSI material).  The components are embedded in both layers.  The wires are also embedded in both layers, but one layer holds only horizontal connections, and the other holds only vertical ones.  If a connection must be made between two components that are not cohorizontal or covertical, then new components (called \textit{vias}) are added to connect horizontal and vertical wires together, resulting in bent wires that alternate between the layers.  However, vias are large compared to wires and their use should be minimized.  In this setting, asking if a graph is a rectangle-visibility graph is the same as asking if a set of components can be embedded so that there is a two-layer routing of their connections that uses no vias.  Our requirement that visibility bands have positive width is motivated by the physical constraint that wires must have some minimum width.  A similar problem arises in printed-circuit board design, as printed-circuit boards naturally have two sides, and connecting wires from one side to the other (the equivalent of making vias) is relatively expensive.}
\end{quote}
In 1995, Dean and Hutchinson began the study of RVGs in their own right.  They focused on bipartite RVGs, and showed that $K_{p,q}$ is an RVG if and only if $p \leq 4$, and that every bipartite RVG with $n$ vertices has at most $4n-12$ edges \cite{Dean95}.  In 1996 Bose, Dean, Hutchinson, and Shermer sought to characterize families of graphs that are RVGs.  They proved that every graph with maximum degree four is an RVG, and every graph that can be decomposed into two caterpillar forests is an RVG, among other results \cite{bose1996rectangle}.  In 1999 Hutchinson, Shermer, and Vince \cite{hutchinson1999representations} proved that every RVG with $n$ vertices has at most $6n-20$ edges, and this bound is tight for $n\geq 8$.  RVGs have since been studied by many other authors \cite{angelini18,biedl16,  Cahit98, GD2014, Dean97,Dean98,dean08,Dean2010,Streinu03}, including generalizations to 3-dimensional boxes \cite{bose99,develin03,Fekete99,gethner11}, rectilinear polygons with more than four edges \cite{digiacomo18, liotta21}, and other variations.

In 1997, Kant, Liotta, Tamassia, and Tollis considered the minimum area, height, and width required to represent a tree as an RVG, as measured by the smallest bounding box containing all of the rectangles in the RV-representation \cite{kant97}.  They obtained asymptotic bounds on the area, width, and height of these representations and found a linear-time algorithm to construct them.  In this paper we consider a similar problem, but seek exact bounds on the area, width, height, and perimeter of an RV-representation of any graph with $n$ vertices.  We say that $\area(G)$, $\perimeter(G)$, $\height(G)$, and $\width(G)$ are the minimum area, perimeter, height, and width, respectively, of the bounding box of any integer rectangle visibility representation of the graph $G$. These are the objects of study in this paper.

In Section~\ref{Definitions} we specify the rectangle visibility graphs we consider and provide definitions and notation needed for the paper.  We finish the section with lemmas we will use in later sections.

In Section~\ref{SeparatingExamples} we show that these four measures of size of a rectangle visibility graph are all distinct, in the sense that there exist two graphs $G_1$ and $G_2$ with $\area(G_1)<\area(G_2)$ but $\perimeter(G_2)<\perimeter(G_1)$, and analogously for all other combinations of these parameters. 

In Section~\ref{SeparatingRepresentations} we show that these measures are not necessarily all attained by the same representation; i.e., there is a graph $G_3$ with two RVG representations $S_1$ and $S_2$ with $\area(G_3)=\area(S_1)<\area(S_2)$ but $\perimeter(G_3)=\perimeter(S_2)<\perimeter(S_1)$, and analogously for all other combinations of these parameters.

In Section \ref{SmallParameters} we characterize the graphs that have the smallest height, width, area, and perimeter among all graphs with $n$ vertices.

In Section~\ref{LargeParameters} we investigate the graphs with largest height, width, area, and perimeter.  We show that, among graphs with $n \leq 6$ vertices, the empty graph $E_n$ has largest area, and for graphs with 7 or 8 vertices, the complete graphs $K_7$ and $K_8$ have larger area than $E_7$ and $E_8$, respectively. Using this, we show that for all $n \geq 7$, the empty graph $E_n$ does not have largest area among all RVGs on $n$ vertices.  The graphs with more than 6 vertices that maximize these parameters are still unknown.

In Section~\ref{Conclusions}, we conclude with a number of open questions.

\section{Basic Definitions and Results} \label{Definitions} 

A rectangle with horizontal and vertical sides whose corners are integer lattice points is said to be an \textit{\textbf{integer rectangle}}.  We consider only integer rectangles for the remainder of the paper. Each rectangle is specified by the two $x$-coordinates and two $y$-coordinates of its corners.  For a set $S$ of rectangles, the smallest rectangle with horizontal and vertical sides containing $S$ is the \textit{\textbf{bounding box}} of $S$.

Suppose $G$ is an RVG with RV-representation $S$ contained in the bounding box $R$, and say $R$ has corners with $x$-coordinates $0$ and $u$ and $y$-coordinates $0$ and $v$, with $u,v\in \Z$.  We can view $R$ as a $u \times v$ grid, with $v$ rows and $u$ columns, and with rectangles in $S$ each contained in a consecutive set of rows and columns.  For example, in the representation shown in Figure~\ref{fig:RVG-example}, rectangle $A$ is contained in rows 3, 4, 5, and 6, and column 1.

We use the convention that lower case letters are vertices of the graph $G$, and the corresponding upper case letters are rectangles in the RV-representation $S$ of $G$; e.g., $a$ is a vertex in $G$, and $A$ is its corresponding rectangle in $S$.  For a given rectangle $A$ in $S$, we denote the $x$-coordinates of its vertical sides by $x_1^A$ and $x_2^A$ with $x_1^A < x_2^A$, and the $y$-coordinates of its horizontal sides by $y_1^A$ and $y_2^A$, with $y_1^A < y_2^A$. In other words, as a Cartesian product of intervals, we have $$A=[x_1^A,x_2^A] \times [y_1^A,y_2^A].$$

We also introduce notation to refer to the set of rectangles in $S$ that are above (respectively, below, to the left of, or to the right of) a given rectangle $A$. Specifically, let the set of rectangles above (north of) $A$ be denoted by
$$ {\mathcal{N}}(A) = \{ X \in S :  y_1^X \geq y_2^A \mbox{ and } (x_1^X,x_2^X) \cap (x_1^A,x_2^A) \not= \emptyset \}.$$  Similarly define ${\mathcal{S}}(A)$, ${\mathcal{W}}(A)$, and ${\mathcal{E}}(A)$ (rectangles south, west, and east of $A$, respectively). For example, in Figure \ref{fig:RVG-example}, $\mathcal{N}(D)=\{B,E\},$ while $\mathcal{E}(A)=\{B, E, D\}$ and $\mathcal{S}(A)=\varnothing.$  Note that $A$ might not see every rectangle in ${\mathcal{N}}(A)$ if there are other rectangles obstructing the view (and similarly for rectangles in the other three sets). 

Let $R$ be the smallest bounding box having horizontal and vertical sides and containing all the rectangles in a set of integer rectangles $S$.  For the remainder of the paper, we turn $R$ so that $\height(R) \leq \width(R)$. 
Given a graph $G$, the \textit{\textbf{area}}, \textit{\textbf{perimeter}}, \textit{\textbf{height}}, and \textit{\textbf{width}} of $G$ are the minimums of the corresponding parameters taken over all bounding boxes of RV-representations of $G$ with height less than or equal to width.

We conclude this section with some preliminary results. 
First we explore the extent to which we can focus on the parameters of connected graphs, and in what ways the values for disconnected graphs are determined or bounded by the parameters of their connected components.

For convenience in stating the next result, we introduce the following notation.
For any positive integers $h$ and $w$, let ${\mathcal{F}}_{h,w}$ denote the (finite) set of graphs that have RV-representations in an $h \times w$ bounding box.


\begin{lemma}\label{NewDisjointLem}
If $G$ is the disjoint union of graphs $H$ and $J$, then:
\begin{enumerate}[label=\rm{(\roman*).}]
\item $\height(G)=\height(H)+\height(J),$
\item $\perimeter(G)=\perimeter(H)+\perimeter(J),$ 
\item $\width(G)=\min \{  \max \{x+b, y+a\} \, | \, H \in {\mathcal{F}}_{x,y}, J \in {\mathcal{F}}_{a,b} \}$,
\item $\area(G)=\min \{ (x+a)(y+b) \, | \, H \in {\mathcal{F}}_{x,y}, J \in {\mathcal{F}}_{a,b} \}$.
\end{enumerate} 
\end{lemma}

\begin{figure}[ht!]
\centering
\includegraphics[width=.7\textwidth]{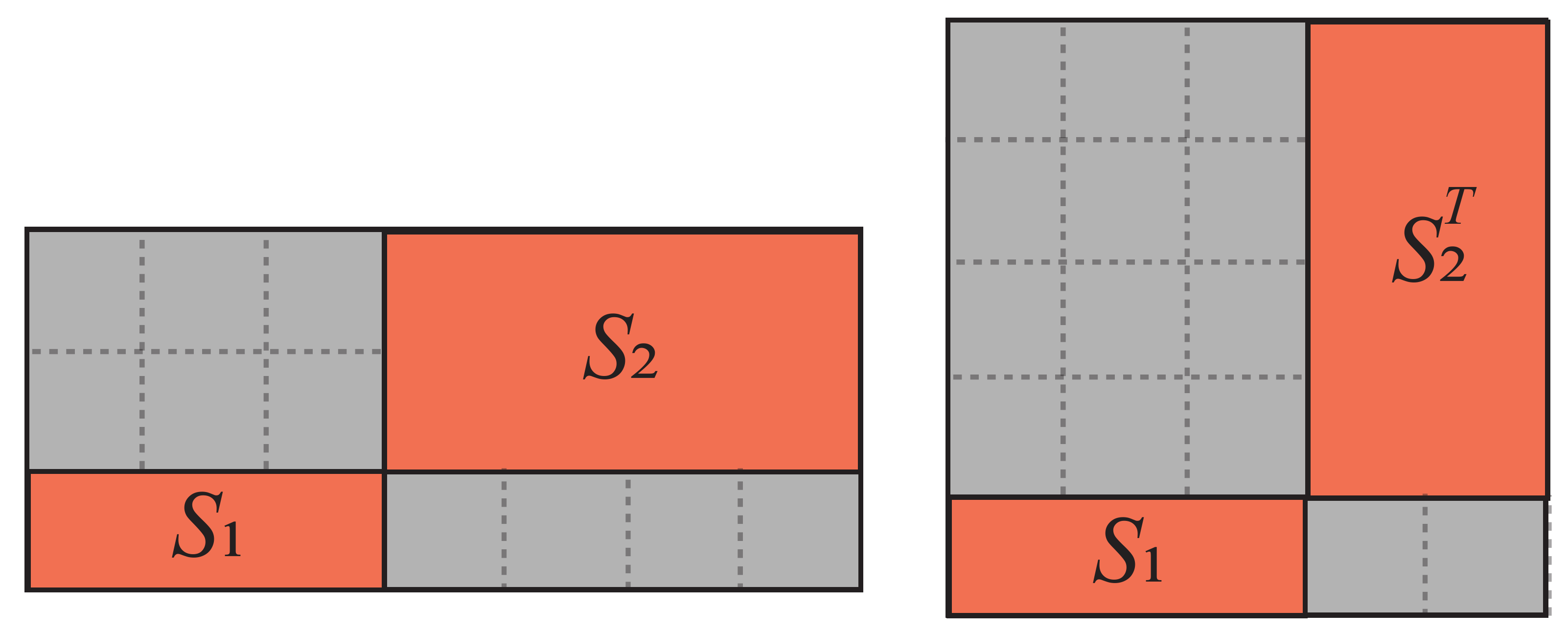}
\caption{Two options for a combined representation of a disjoint union of two RVGs}
\label{Glue} 
\end{figure}

\begin{proof}
Suppose $G$ is the disjoint union of graphs $H$ and $J$.  Given any RV-representations $S_1$ and $S_2$ of $H$ and $J$, we construct two RV-representations of $G$. As indicated in Figure~\ref{Glue}, we identify the upper right corner of $S_1$ with the lower left corner of either $S_2$ or  $S_2^T$, where $S_2^T$ denotes the RV-representation of $J$ formed by transposing $S_2$ across its main (top left to lower right) diagonal.  If $S_1$ is $x \times y$ and $S_2$ is $a \times b$, then it follows that
$ G \in {\mathcal{F}}_{x+a,y+b} \cap {\mathcal{F}}_{x+b,y+a}.$
This implies each of the expressions in (i)-(iv) are upper bounds for the parameters of $G$.

To prove these expressions are also lower bounds, note that 
 any RV-representation $S$ of $G$ must have the rectangles corresponding to vertices of $H$ in separate rows and columns from the rectangles corresponding to vertices of $J$.  If $S$ has height smaller than $\height(H)+\height(J)$, then either the rectangles in $S$ corresponding to vertices of $H$ must form an RV-representation of $H$ with height less than $\height(H)$ or the rectangles in $S$ corresponding to vertices of $J$ must form an RV-representation of $J$ with height less than $\height(J)$.  Therefore no such representation is possible.  Similar arguments apply to the other parameters.
\end{proof}

\begin{remark}
The following example illustrates a subtlety captured by the formula in Lemma~\ref{NewDisjointLem}(iv).  Consider the graphs $P_4+K_{1,4}$ and $K_{1,4}+K_{1,4}$.  We have $\area(P_4+K_{1,4})=27$, obtained from a $1 \times 4$ RV-representation of $P_4$ and a $2 \times 5$ RV-representation of $K_{1,4}$.  But we also have $\area(K_{1,4}+K_{1,4})=36$, obtained from two copies of a $3 \times 3$ RV-representation of $K_{1,4}$ (both representations of $K_{1,4}$ are given in Figure~\ref{fig:width-examples2a}).  So the minimum area of $H + K_{1,4}$ uses different representations of $K_{1,4}$ depending on $H$.
\end{remark}

\begin{corollary}\label{NewDisjointCor}
If $G$ is the disjoint union of graphs $H$ and $J$, then:
\begin{enumerate}[label=\rm{(\roman*).}]
\item $\width(G) \leq \width(H) + \width(J)$,
\item $\area(G) \leq (\width(H) + \width(J))^2$.
\end{enumerate} 
\end{corollary}

\begin{proof}
These both follow immediately from the 
construction in the proof of Lemma~\ref{NewDisjointLem} shown in Figure~\ref{Glue}. 
\end{proof}

\begin{lemma}
Suppose $G$ is a graph with RV-representation $S$ with bounding box $R$.  If $\height(G)=\height(R)$ and $\width(G)=\width(R)$, then $\area(G)=\area(R)$ and $\perimeter(G)=\perimeter(R)$.  In other words, if $S$ realizes the height and width of $G$, then $S$ also realizes the area and perimeter of $G$.  
\label{lemma-same-height-width}
\end{lemma}

\begin{proof}
Any representation with smaller area or perimeter must have smaller height or smaller width, which is impossible by hypothesis.
\end{proof}

Later (in Table~\ref{table-separating2}), we will see that the hypotheses of Lemma~\ref{lemma-same-height-width} are necessary.  In particular, $G_4$ has two different representations for minimizing area and perimeter.


\section{Height, Width, Area, and Perimeter induce distinct orderings of RVGs} \label{SeparatingExamples} 

\begin{figure}[ht!]
    \centering \includegraphics[width=.9\textwidth]{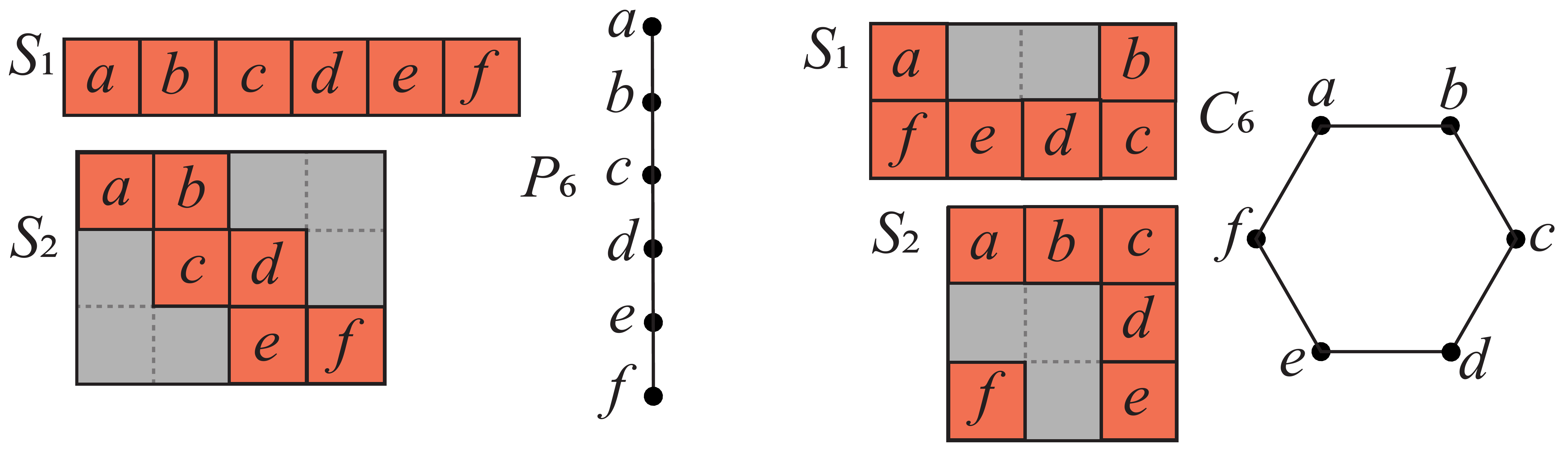}
    \caption{Graphs $P_6$ and $C_6$ show that height and area order graphs differently than  width and perimeter.}
    \label{fig:width-examplesAA}
\end{figure}

\begin{figure}[ht!]
    \includegraphics[width=.8\textwidth]{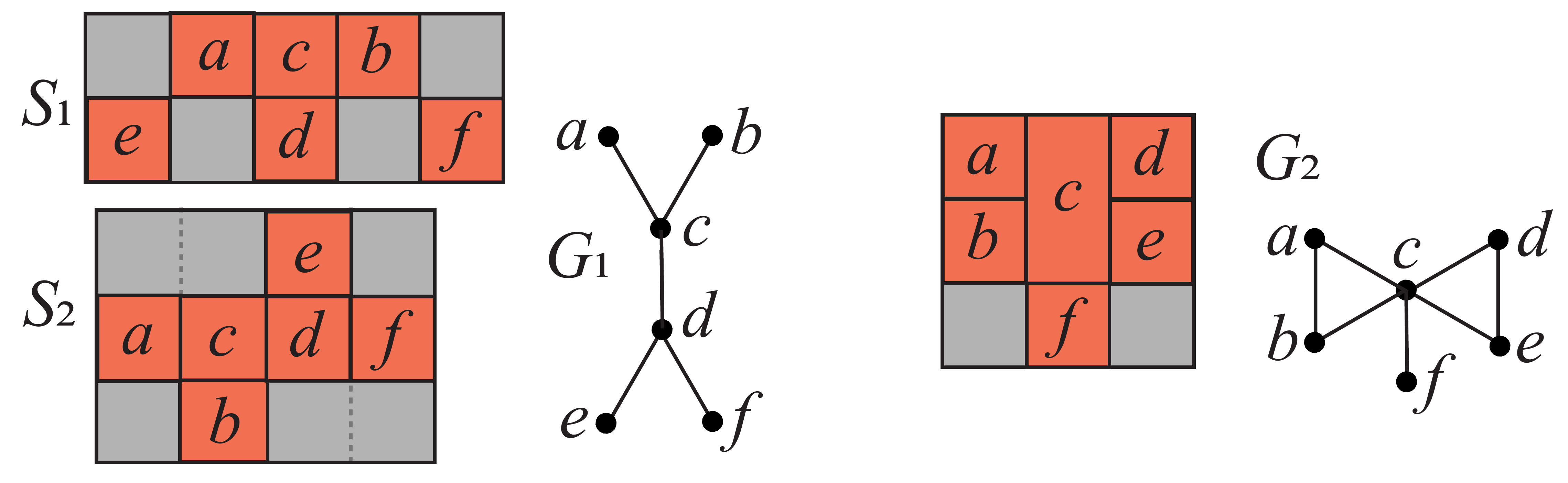}
    \caption{Graphs $G_1$ and $G_2$ show that height orders graphs differently than area.}
    \label{fig:width-examplesBB}
\end{figure}

\begin{figure}[ht!]
    \centering
    \includegraphics[width=.4\textwidth]{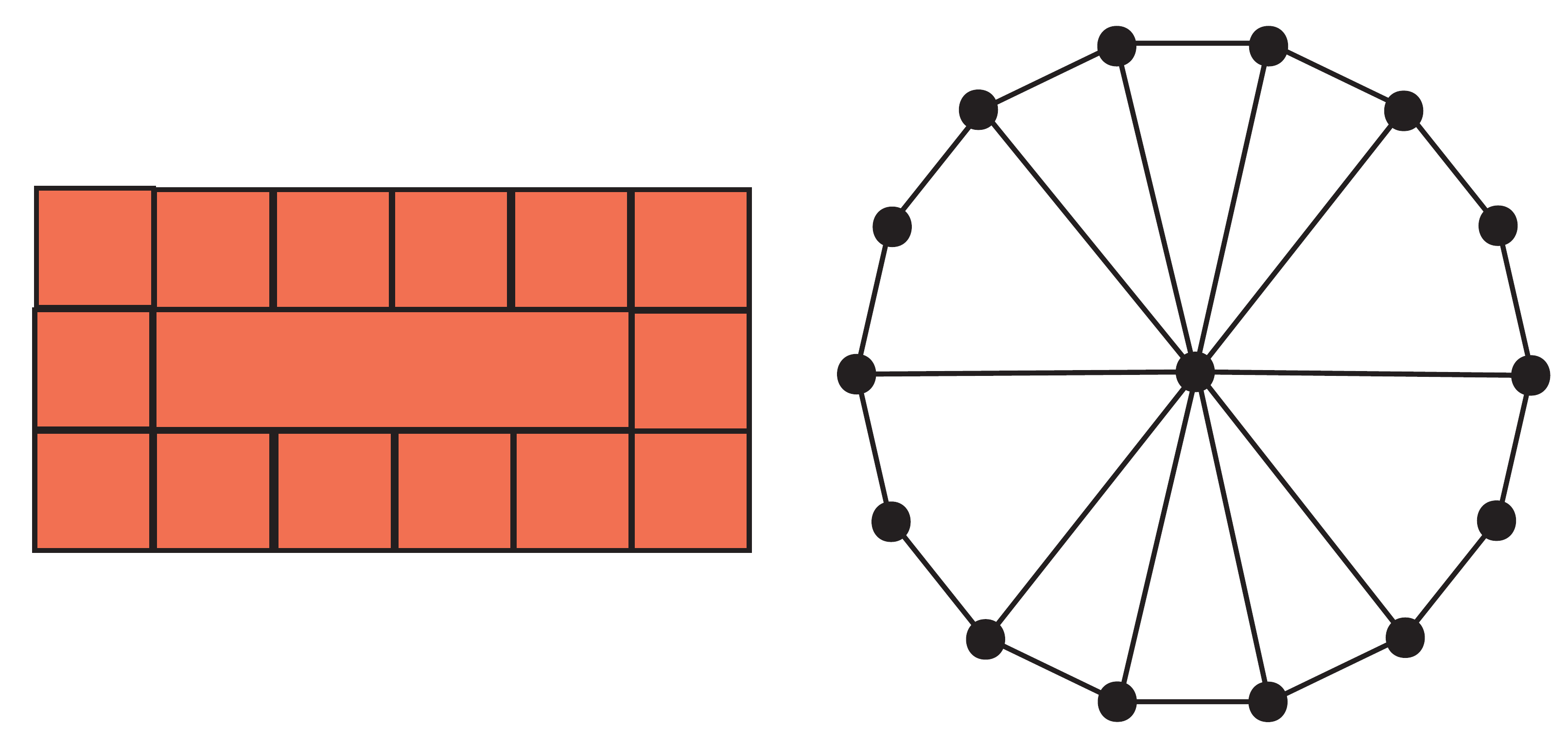} \hspace{1em}
    \includegraphics[width=.4\textwidth]{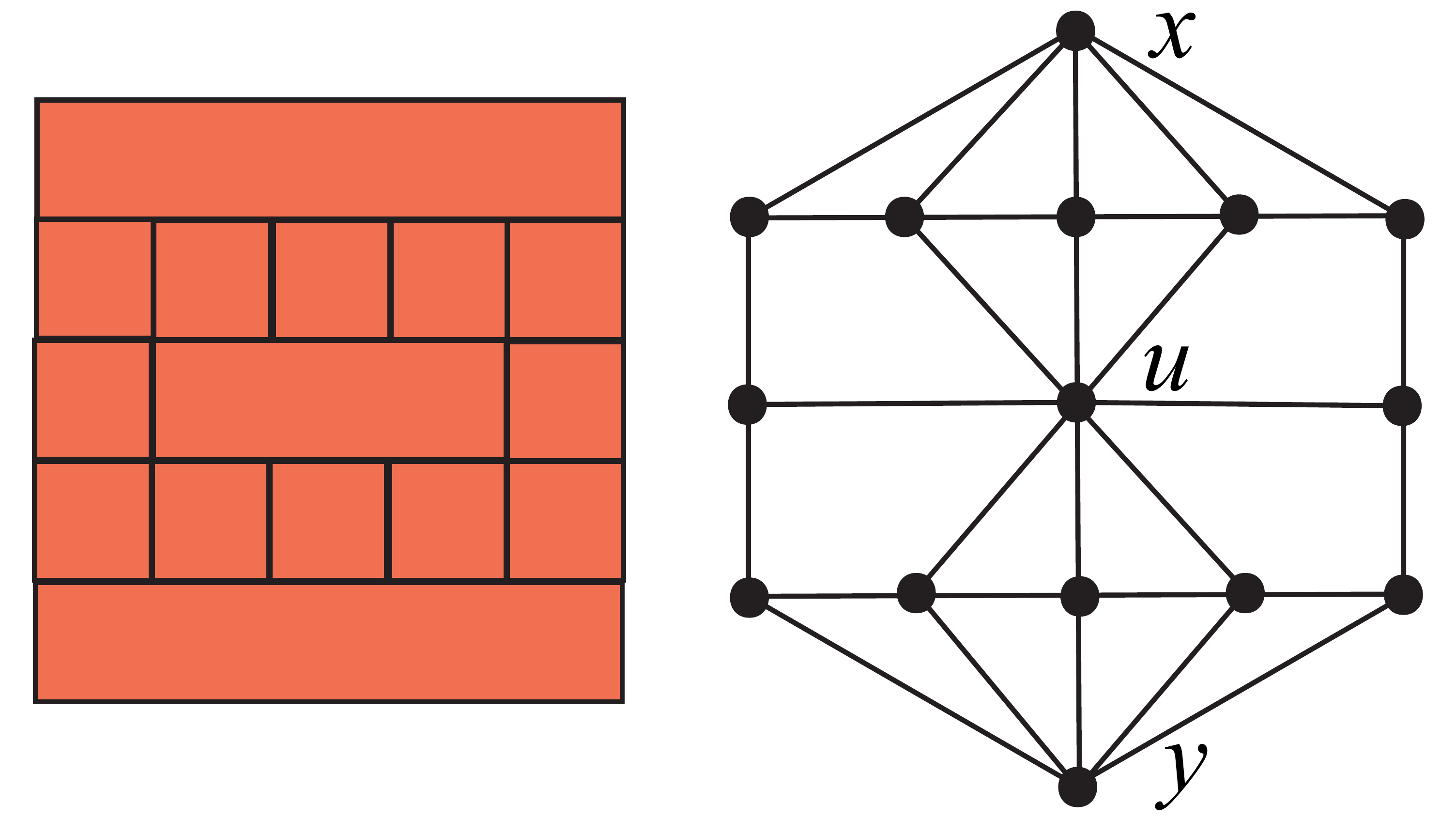}
    \caption{Graphs $G_3$ and $G_4$ show that width orders graphs differently than perimeter.}
    \label{fig:width-examplesCC}
\end{figure}

In this section, we consider the various notions of height, width, perimeter, and area of RVGs. We show that these parameters represent independent measures of RVGs, in the sense that they do not always give identical orderings of the sets of graphs on a given number of vertices. 
Examples to illustrate these results are summarized in Tables~\ref{table-separating1}  and \ref{table-separating-new}. Graphs $G_1$, $G_2$, $G_3$, and $G_4$ in these tables are shown in Figures~\ref{fig:width-examplesBB} and~\ref{fig:width-examplesCC}.  For each pair of parameters, there exists a pair of graphs with an equal number of vertices that are oppositely ordered by those parameters. We have verified the height, width, area, and perimeter of all connected graphs with 6 or fewer vertices by computer search   \cite{WebList} and the claims regarding $P_6$, $C_6$, $G_1$, and $G_2$ in Tables~\ref{table-separating1}  and \ref{table-separating-new} follow easily, see Figures~\ref{fig:width-examplesAA} and \ref{fig:width-examplesBB}.
The claims regarding $G_3$ and $G_4$, each with 15 vertices, are proved in Theorems~\ref{width-theorem} and \ref{width-theorem2}.


\begin{table}[ht!]
\begin{tabular}{c|c|cccc} 
 Graph & Vertices & Height & Width & Area & Perimeter \\ \hline \hline
  $P_6$  & 6 & {\bf 1} & 4 & {\bf 6} & 14 \\ 
  $C_6$  & 6 & 2 & {\bf 3} & 8 & {\bf 12} \\ \hline
  $G_1$ & 6 & {\bf 2} & -- & 10 & -- \\
  $G_2$ & 6 & 3 & -- & {\bf 9} & -- \\ \hline
 $G_3$ &  15 & -- & 6 & -- & {\bf 18} \\
 $G_4$ & 15 & -- & {\bf 5} & -- & 20 \\   \hline
\end{tabular} 
\smallskip
\caption{Separating examples for height, width, area, and perimeter.}
\label{table-separating1}
\end{table}

\begin{table}[ht!]
\begin{tabular}{c|ccc}
 & Perimeter & Height & Width \\ \hline \hline
  Area & $P_6,C_6$ & $G_1,G_2$ & $P_6,C_6$ \\ 
  Perimeter & -- & $P_6,C_6$ & $G_3,G_4$ \\ 
  Height & -- & -- & $P_6,C_6$ \\
  \hline
\end{tabular} 
\smallskip
\caption{Graph pairs that are oppositely ordered by each of the parameter pairs.}
\label{table-separating-new}
\end{table}

\begin{theorem}
The graph $G_3$ shown on the left in Figure~\ref{fig:width-examplesCC} has perimeter 18 and width 6. \label{G1-proof}
\label{width-theorem}
\end{theorem}

\begin{proof}
\noindent To find the perimeter of $G_3$, we first consider its area. Let $v$ be the vertex of degree 10 in $G_3$. In any RV-representation of $G_3$, the corresponding rectangle $V$ must have perimeter at least 10, so its area is at least 4. Together with the 14 other vertices, we see $\area(G_3)\ge 18.$  Now suppose $G_3$ can be represented in a bounding box of height $h$, width $w$, and perimeter $p=2h+2w$. Such a rectangle has maximum area when $h=w=p/4$, so $\area(R) \leq p^2/16$. But $\area(R) \geq 18$, so $p \geq \sqrt{18 \cdot 16} > 16.$  Since $p$ must be even,  $\perimeter(G_3)=18$ by Figure \ref{fig:width-examplesCC}.

Next we consider $\width(G_3)$.  If $\width(G_3)<6$, then $G_3$ can be represented in a $5\times 5$ box $R$.   In a $5 \times 5$ box with 14 other vertices, $V$ has area at most 11, and hence $V$ must be $1 \times 4$, $1 \times 5$, $2 \times 3$, $2 \times 4$, $2 \times 5$, or $3 \times 3$.  We rule out each possibility below.

\begin{figure}[ht!]
    \centering
    \includegraphics[width=.9\textwidth]{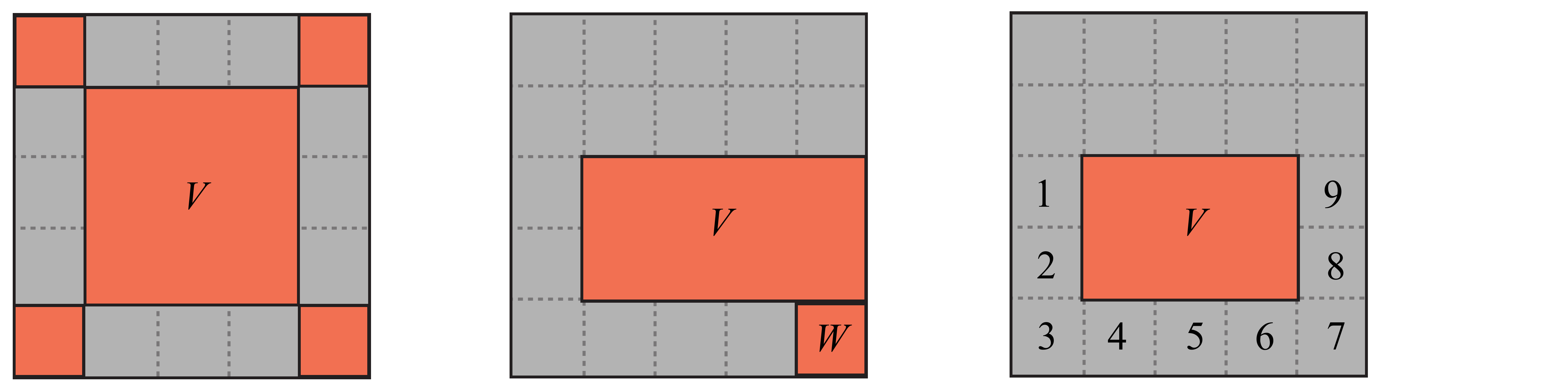}
    \caption{Diagrams illustrating several of the cases in the proof of Theorem~\ref{G1-proof}.}
    \label{fig:width-proof}
\end{figure} 
 
  $\bullet \;$ If  $V$ is $1 \times 4$, at least one unit of the perimeter of $V$ is on the boundary of $R$, and therefore $V$ cannot represent a degree-10 vertex. 

  $\bullet \;$ If $V$ is $1 \times 5$ or $2 \times 5$, it must touch opposite sides of $R$; thus $V$ must see 5 rectangles on each of its other two sides in order to have degree 10. But then $v$ is a cut vertex, and $G_3$ has none. 

  $\bullet \;$ If  $V$ is $3\times 3$, it cannot have an edge on the boundary of $R$ and represent a degree-10 vertex. In this case, $V$ occupies the middle of $R$ as shown on the left in Figure~\ref{fig:width-proof}. The four vertices not adjacent to $v$ must be represented by $1\times 1$ squares in the corners of $R$.  Among these four, two (disjoint) pairs have a unique common neighbor in $G_3$.  Locating the rectangles for these common neighbors in their respective  $3 \times 1$ blocks of $R$, we see there now remain only 6 locations for the other 8 vertices.   

  $\bullet \;$ If $V$ is $2\times 4$, it must (by symmetry) appear as in the middle of Figure~\ref{fig:width-proof}. To have degree 10, $V$ must see a $1\times 1$ square at $W$. But $v$ has no neighbor of degree less than $3$.   

  $\bullet \;$ If $V$ is $2\times 3$, it must (by symmetry) appear as on the right of Figure~\ref{fig:width-proof}.  To have degree 10, $V$ must see distinct rectangles on each unit of its perimeter.  Numbering the locations in $R$ as in the right of Figure~\ref{fig:width-proof}, if location 3 is empty, location 2 must contain a $1\times 1$ square, but $v$ has no neighbor of degree less than $3$. If a $2 \times 1$ rectangle $X$ covers location 3, the rectangles in locations 1 through 6 form a path $P_5$ in the neighborhood of $V$, but $G_3$ has no such subgraph.  Thus, locations 3 and (by symmetry) 7 contain $1 \times 1$ squares. But $G_3 -v$ has no vertices of degree 2 that are at distance 4.
 \end{proof}

\begin{theorem}
The graph $G_4$ shown on the right in Figure~\ref{fig:width-examplesCC} has perimeter 20 and width 5. \label{G2-proof}
\label{width-theorem2}
\end{theorem}
               
\begin{proof} 
We use the labeling in Figure~\ref{fig:width-examplesCC}.  We begin with perimeter.  Note that since $G_4$ has 15 vertices, every RV-representation of $G_4$ has area at least 15.  Then if $\perimeter(G_4)<20$, it must have a representation that fits in a $3 \times 6$ or $4 \times 5$ bounding box. But a $3 \times 6$ box has area 18, so the 3 vertices of degree more than four in $G_4$ must be represented with $1 \times 2$ rectangles. But then all the vertices of $G_4$ must have degree $\leq 6$, a contradiction.

Thus $G_4$ must be representable in a $4 \times 5$ box. The two vertices of degree 5 in $G_4$ must each be represented with a rectangle of area at least 2. In a $4 \times 5$ box $R$, $U$ then has area 3 or 4 and cannot lie in a corner. All vertices of $G_4$ have degree at least 3, so no corner of $R$ has a $1 \times 1$ square. Therefore, each corner of $R$ is either empty or occupied by a rectangle of area at least 2. With 15 vertices, this exceeds the total area of 20 in $R$: each vertex contributes at least one, the four corners require at least one additional unit of area each, and $U$ needs two or three additional units of area not lying in a corner. Now $\perimeter(G_4)=20$ by Figure~\ref{fig:width-examplesCC}.

We next consider width. If $\width(G_4)<5$, then $G_4$ can be represented in a $4 \times 4$ box. But then it also has a $4 \times 5$ representation, which we just proved impossible. By Figure~\ref{fig:width-examplesCC}, we see that $\width(G_4) =5$. 
\end{proof}


\section{Minimizing Height, Width, Area, and Perimeter can require distinct representations of an RVG} \label{area} \label{SeparatingRepresentations} 

\begin{figure}[ht!]
    \centering
   \includegraphics[width=\textwidth]{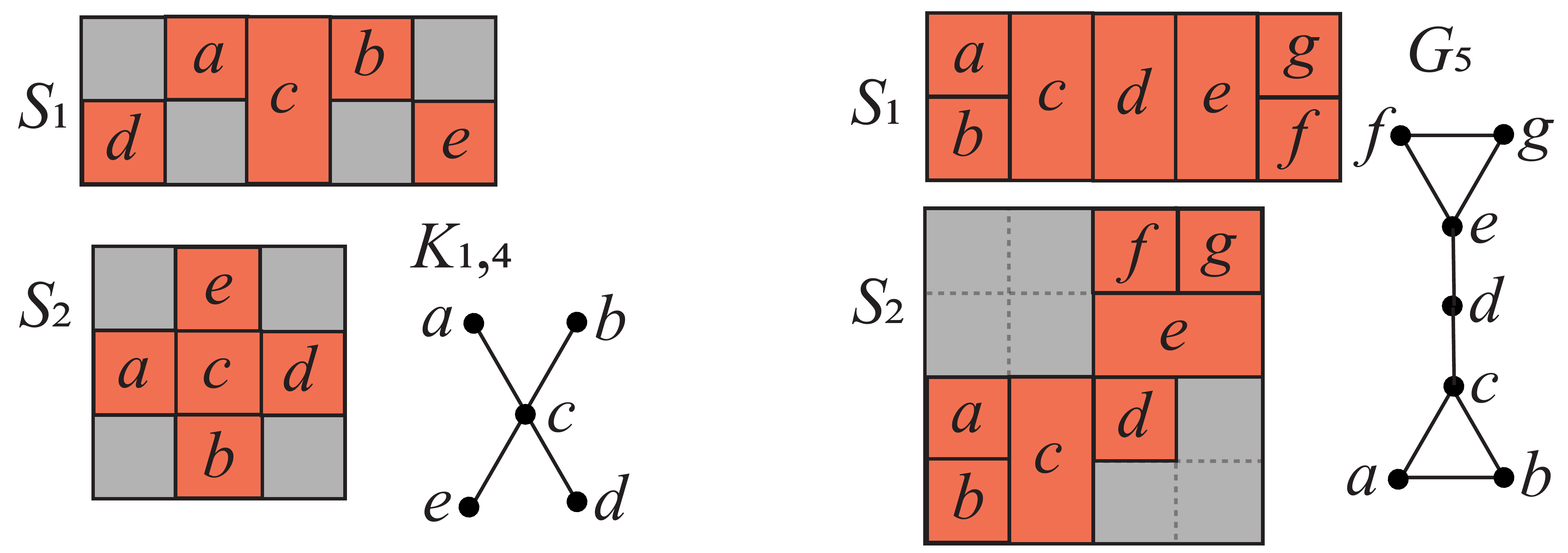}
    \caption{Graphs $K_{1,4}$ and $G_5$ are shown.  For $K_{1,4}$, one representation ($S_1$) minimizes height and another ($S_2$) minimizes area, perimeter, and width.  For $G_5$, one representation ($S_1$) minimizes height, area, and perimeter, and another  ($S_2$) minimizes width.}
    \label{fig:width-examples2a}
\end{figure}

\begin{figure}[ht!]
    \centering
   \includegraphics[width=.6\textwidth]{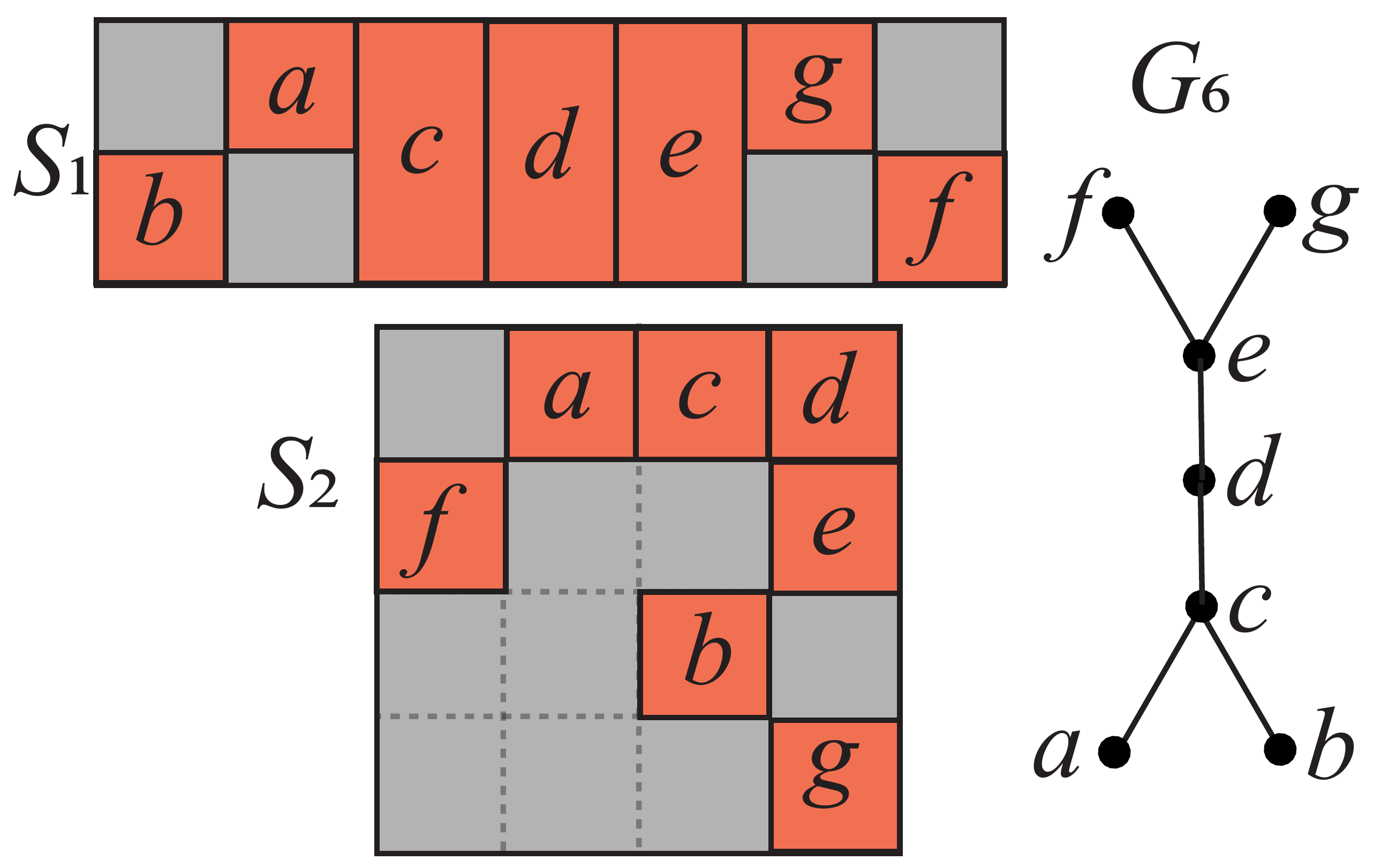}
    \caption{Graph $G_6$ with distinct representations minimizing area and perimeter.}
    \label{fig:width-examples2b}
\end{figure}

In this section, we further explore our four parameters and we observe that, even for a single RVG, it is possible that the set of representations minimizing one of them may be disjoint from the set of representations minimizing another.  
Examples to illustrate these results are summarized in Tables~\ref{table-separating2} and \ref{table-separating2b}. For each pair of parameters, there exists a graph that requires distinct representations to separately minimize each parameter in that pair. Specifically, the star $K_{1,4}$ has representations minimizing height that are distinct from those minimizing width, area, and perimeter; the graph $G_5$ shown in Figure~\ref{fig:width-examples2a} has representations minimizing width that are distinct from those minimizing the other parameters; and the graph $G_6$ shown in Figure~\ref{fig:width-examples2b} has distinct representations minimizing area and perimeter. Since $K_{1,4}$ has fewer than 7 vertices, the claims regarding it are verified by our computer search \cite{WebList}. The claims regarding $G_5$ and $G_6$, each with 7 vertices, are proved in Theorems~\ref{G5-theorem} and \ref{G6-theorem} below.  

\begin{table}[ht!]
\begin{tabular}{c|c|cccc} 
Graph (representation) & Vertices & Height & Width  & Area & Perimeter\\ \hline \hline 
  $K_{1,4}$ ($S_1$) & 5 & {\bf 2} & 5  & 10 & 14\\
  $K_{1,4}$ ($S_2$) & 5 & 3 & {\bf 3}  & {\bf 9} & {\bf 12}\\ \hline
  $G_5$ ($S_1$) & 7 & {\bf 2} & 5  & {\bf 10} & {\bf 14}\\
  $G_5$ ($S_2$) & 7 & 4 & {\bf 4}  & 16 & 16\\ \hline
  $G_6$ ($S_1$) & 7 & {\bf 2} & 7  & {\bf 14} & 18\\
  $G_6$ ($S_2$) & 7 & 4 & {\bf 4}  & 16 & {\bf 16}\\  \hline
\end{tabular} 
\smallskip
\caption{Graphs representations used to minimize height, width, area, and perimeter.}
\label{table-separating2}
\end{table} 

\begin{table}[ht!]
\begin{tabular}{c|ccc}
 & Perimeter & Height & Width \\ \hline \hline
  Area & $G_6$ & $K_{1,4}$ & $G_5$ \\ 
  Perimeter & -- & $K_{1,4}$ & $G_5$ \\ 
  Height & --& --& $G_5$ \\
  \hline
\end{tabular} 
\smallskip
\caption{Graphs requiring different representations to minimize area, perimeter, height, and width.}
\label{table-separating2b}
\end{table}

\begin{theorem}
The graph $G_5$ shown in Figure~\ref{fig:width-examples2a} has height 2, width 4, area 10, and perimeter 14.  Furthermore, minimizing width requires a different RV-representation from height, area, or perimeter. \label{G5-theorem}
\end{theorem}

\begin{proof}
 The representation $S_1$ in Figure \ref{fig:width-examples2a} has height 2, width 5, area 10, and perimeter 14. The representation $S_2$ has height 4, width 4, area 16, and perimeter 16. 
 
We claim $G_5$ has no $3 \times 4$ representation. Each of the 3-cycles in $G_5$ must have at least 2 rows that occupy at least 2 units of area. This implies that some row in a $3 \times 4$ box has 2 units of each 3-cycle, which therefore must see each other. Since $G_5$ has no edges joining these 3-cycles, this is impossible.

It follows that $G_5$ has no $3 \times 3$ or $2 \times 4$ representation. Since $G_5$ is not a path, $\height(G_5)\geq 2$ and these facts together imply that $G_5$ has height 2, width 4, area 10, and perimeter 14.  
Since the only representations of $G_5$ with width 4 are $4 \times 4$, any representation minimizing width does not minimize height, area, or perimeter.
\end{proof}

\begin{theorem}  
The graph $G_6$ shown in Figure~\ref{fig:width-examples2b} has perimeter 16 and area 14.  Furthermore, minimizing perimeter requires a
different RV-representation than area. \label{G6-theorem}
\end{theorem}

\begin{proof}
 The first representation in Figure \ref{fig:width-examples2b} has perimeter 18 and area 14. The second has perimeter 16 and area 16. 
 
 The 4 vertices of degree 1  each require at least 3 units of length on the perimeter with free lines of sight. The degree 2 vertex requires at least 2 units, and the 2 degree 3 vertices each require at least 1 unit. So $\perimeter(G_6)=16$.  
 
 If $\area(G_6)<14,$ then $G_6$ can be represented in a rectangle of area 7, 8, ..., or 13. Since $G_6$ is not a path, prime areas are not possible, so $G_6$ must have height $>1$. Since $\perimeter(G_6)=16$, the only possible bounding box must be $2 \times 6$. We claim this is impossible. Since $G_6$ has 7 vertices, a box with 6 columns would force some column to contain portions of 2 rectangles. Neither of these rectangles can represent a vertex of degree 1, or else the remaining graph must be a path, so they must have degrees 2 and 3, which implies that $D$ has an empty horizontal line of sight. But now, taken together, the 5 vertices of degree 1 or 2 require 5 empty horizontal lines of sight, while the bounding box only has height 2, with at most 4 such lines.   
 
 The only rectangle with height $>1$ and area 14 is $2 \times 7$, so the perimeter and area must be achieved with distinct representations.
\end{proof}


\section{Graphs with small area, perimeter, height, and width} \label{SmallParameters}  

In this section, we address the question of which graphs on a given number of vertices minimize each parameter. 

For any real number $x$, we use $\lceil x \rceil$ and $\lfloor x \rfloor$ to denote the integer ceiling and floor of $x$, respectively. We use $[[ x ]] = \lfloor x + 1/2 \rfloor$ to denote $x$ rounded to the nearest integer. Recall that, for any positive integers $h$ and $w$, we let ${\mathcal{F}}_{h,w}$ denote the (finite) set of graphs that have RV-representations in an $h \times w$ bounding box.

\begin{theorem}\label{lowerbds}
Let $G$ be any graph with $n$ vertices and suppose $G$ has an RV-representation. Then the following hold.
\begin{enumerate}[label=\rm{(\roman*).}]
    \item The height of $G$ satisfies
    $$\height(G) \geq 1.$$
     Equality holds if and only if $G \cong P_n$, the path on $n$ vertices.
    \item The area of $G$ satisfies
    $$\area(G) \geq n.$$  
     Equality holds if and only if $G \cong P_h \openbox P_w$, for some positive integers $h$ and $w$, where $n=h \cdot w$. 
    \item The width of $G$ satisfies
    $$\width(G) \geq \lceil \sqrt{n} \rceil.$$ 
    Equality holds if and only if $G \in {\mathcal{F}}_{w,w}$ where $w = \lceil \sqrt{n} \rceil$.
    \item The perimeter of $G$ satisfies
    $$\perimeter(G) \geq 2 \cdot [[\sqrt{n}]] + 2 \cdot \lceil \sqrt{n} \rceil.$$ 
    Equality holds if and only if $G \in {\mathcal{F}}_{h,w}$ for some positive integers $h$ and $w$, where $h+w = [[\sqrt{n}]] + \lceil \sqrt{n} \rceil$ and $hw \geq n$.
\end{enumerate} 
\end{theorem}

\begin{proof}
 To see (i) and (ii), note that every RV-representation of a graph must have height at least 1 and area at least $n$, since each vertex requires at least a $1\times 1$ rectangle. For $G$ to have height exactly 1, every rectangle in such a representation must be on the same row in the representation, so $G$ is a path.  For $G$ to have area exactly $n$,  every rectangle in such a representation must be $1 \times 1$, with no empty space in the bounding box $R$. Thus $G$ is the grid $P_h \openbox P_w$, where $h$ is the height of $R$ and $w$ is the width of $R$.  
 
 We show (iii) by way of contradiction.  Suppose $G$ can be represented in an $a \times b$ bounding box $R$, where $a \leq b < \lceil \sqrt{n} \rceil$. Then $ b \leq \lceil \sqrt{n} \rceil -1$, so $b <  \sqrt{n}$. But now the area of $R$ is $ab <  n$, which is impossible since $G$ has $n$ vertices.  It follows  that $\width(G) = \lceil \sqrt{n} \rceil$ if and only if $G \in {\mathcal{F}}_{w,w}$ where $w = \lceil \sqrt{n} \rceil$. 
 
 We also show (iv) by way of contradiction.  Suppose $G$ can be represented in an $a \times b$ bounding box $R$, where $a + b < [[\sqrt{n}]] + \lceil \sqrt{n} \rceil$. By (ii), we know $ab \geq n$. We consider cases for whether $\sqrt{n}$ rounds up or down. In each case, we find that
 $$ (a+b)^2 < 4n \leq 4ab,$$
 which implies that $(a-b)^2 <0,$ a contradiction.  
 It follows that $\perimeter(G) = 2 \cdot [[\sqrt{n}]] + 2 \cdot \lceil \sqrt{n} \rceil$ if and only if $G \in {\mathcal{F}}_{h,w}$ for some positive integers $h$ and $w$, where $h+w = [[\sqrt{n}]] + \lceil \sqrt{n} \rceil$ and $hw \geq n$. 
\end{proof}

\begin{remark}
The condition for equality in Theorem~\ref{lowerbds}(iv) restricts the bounding box to be very nearly square. Specifically, we can say the following. 
Let $k$ denote the integer such that $k^2 < n \leq (k+1)^2$.

If $k^2 < n \leq k(k+1)$, then equality holds in Theorem~\ref{lowerbds}(iv) if and only if $G \in {\mathcal{F}}_{k-t,k+1+t}$ for some integer $t$ where $$0 \leq t \leq \sqrt{(k+{\textstyle\frac{1}{2}})^2-n}-{\textstyle\frac{1}{2}}.$$

If $k(k+1) < n \leq (k+1)^2$, then equality holds in Theorem~\ref{lowerbds}(iv) if and only if $G \in {\mathcal{F}}_{k+1-t,k+1+t}$ for some integer $t$ where $$0 \leq t \leq \sqrt{(k+1)^2-n}.$$
\end{remark}

For example, any RVG with $n=70$ vertices must have perimeter at least 34, with equality only when $G$ has a $7 \times 10$ or $8 \times 9$ representation.
Similarly, any RVG with $n=120$ vertices must have perimeter at least 44, with equality only when $G$ has a $10 \times 12$ or $11 \times 11$ representation. 


\section{Graphs with large area, perimeter, height, and width} \label{CompleteArea} \label{LargeParameters} 

In this section we turn to the question of which graphs on a given number of vertices maximize our four parameters.  

Recall that the \textit{\textbf{empty graph}} $E_n$ is the graph with $n$ vertices and no edges. Among small graphs (at most 6 vertices), the empty graphs maximize each of the four parameters. When the number of vertices is larger than 7, however, we will see that the empty graph no longer reigns supreme. Our proof is constructive, as we will provide specific graphs that we will prove are larger than $E_n$ in each parameter. But our results here leave open, perhaps for future work, the more difficult question of which graphs on $n$ vertices actually achieve the maximum values for the four parameters.

We begin with the small graphs.

\begin{theorem}
For $1 \leq n \leq 6$, among all graphs with $n$ vertices, the empty graph $E_n$ has largest height, width, area, and perimeter.
\end{theorem}

\begin{figure}
\centering
\includegraphics[width=.7\textwidth]{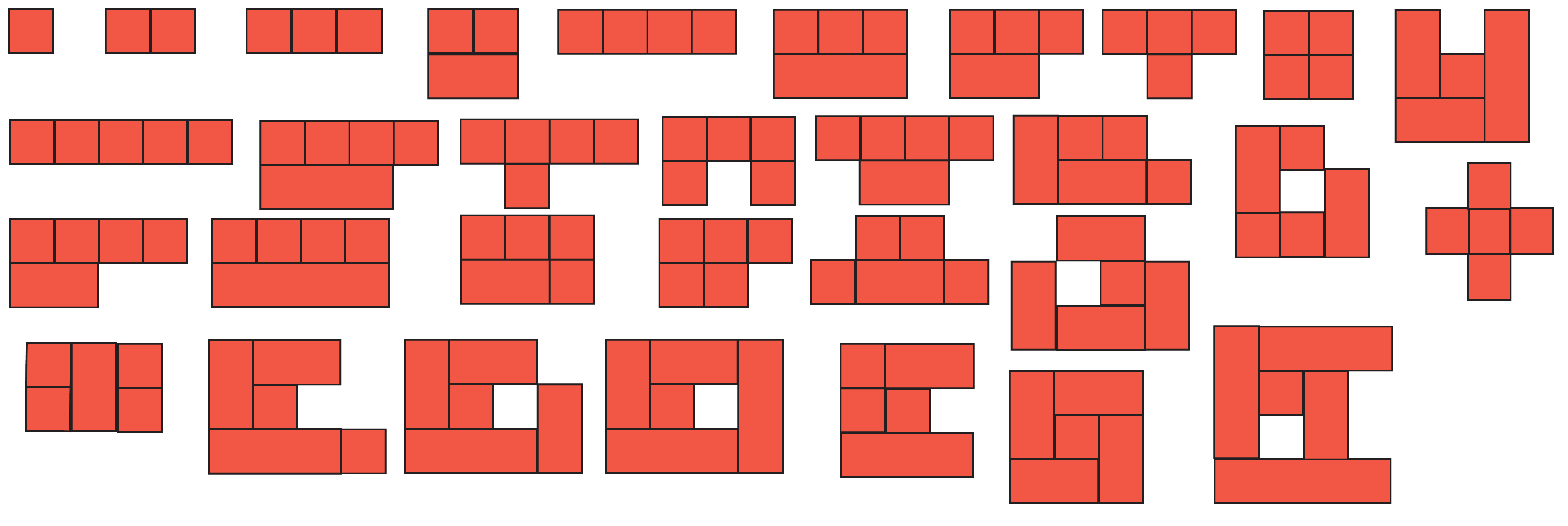}
\caption{RV-representations of all connected graphs with between 1 and 5 vertices.}
\label{fig:small representations}
\end{figure}

\begin{figure}[ht!]
\centering
\includegraphics[width=.9\textwidth]{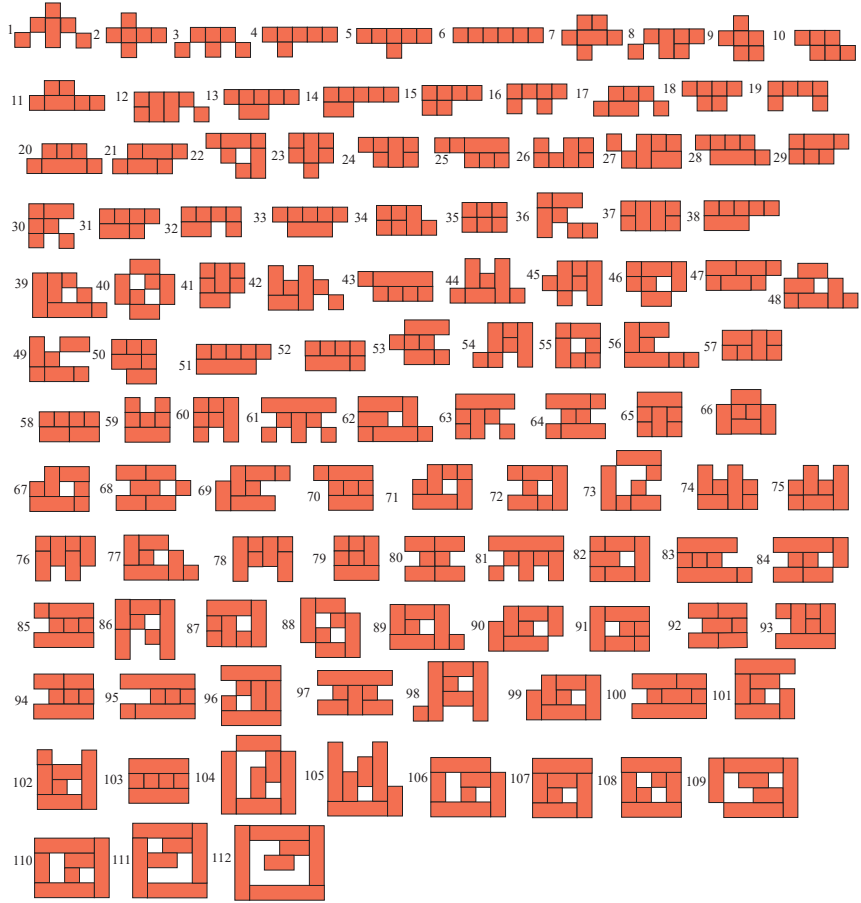}
\caption{RV-representations of all connected graphs with 6 vertices, using the labeling from~\cite{small-graphs-website}.  These representations have smallest area, by computer search.}
\label{fig:small representations6}
\end{figure}

\begin{proof}
By Lemma~\ref{NewDisjointLem}, the empty graph $E_n$ has height $n$, width $n$, area $n^2$, and perimeter $4n$.  Figures~\ref{fig:small representations} and \ref{fig:small representations6} show RV-representations of all connected graphs with at most 6 vertices.  These figures show that no other connected graphs with $2 \leq n \leq 6$ vertices exceed any of these values.  For a disconnected graph $G$, Lemma~\ref{NewDisjointLem} implies that, as long as $G$ has at least one component with more than one vertex, we can combine the representations of the components as in Figure~\ref{Glue} to obtain height less than $n$, area less than $n^2$, and perimeter less than $4n$. Because $K_2$ has width 2, the graph $K_2 + E_{n-2}$ has width $n$, but no graph has larger width than $E_n$ for $n \leq 6$.
\end{proof}

Our next results focus heavily on the RV-representations of the complete graph $K_n$.  Let $S$ be any set of rectangles representing $K_n$ and let $R = [0,u] \times [0,v]$ denote the smallest bounding box containing them. 
%
%
We define the set ${\mathcal{T}}_S$ of {\bf top rectangles of $S$} as follows:
$$
{\mathcal{T}}_S = \{ X \in S :  y_1^X \geq y_1^Y \mbox{ for all } Y \in S\}.$$
We now prove that for $K_n$, ${\mathcal{T}}_S$ contains a single rectangle when $n\geq 6$.

\medskip

\begin{figure}
    \centering
    \includegraphics[height=2in]{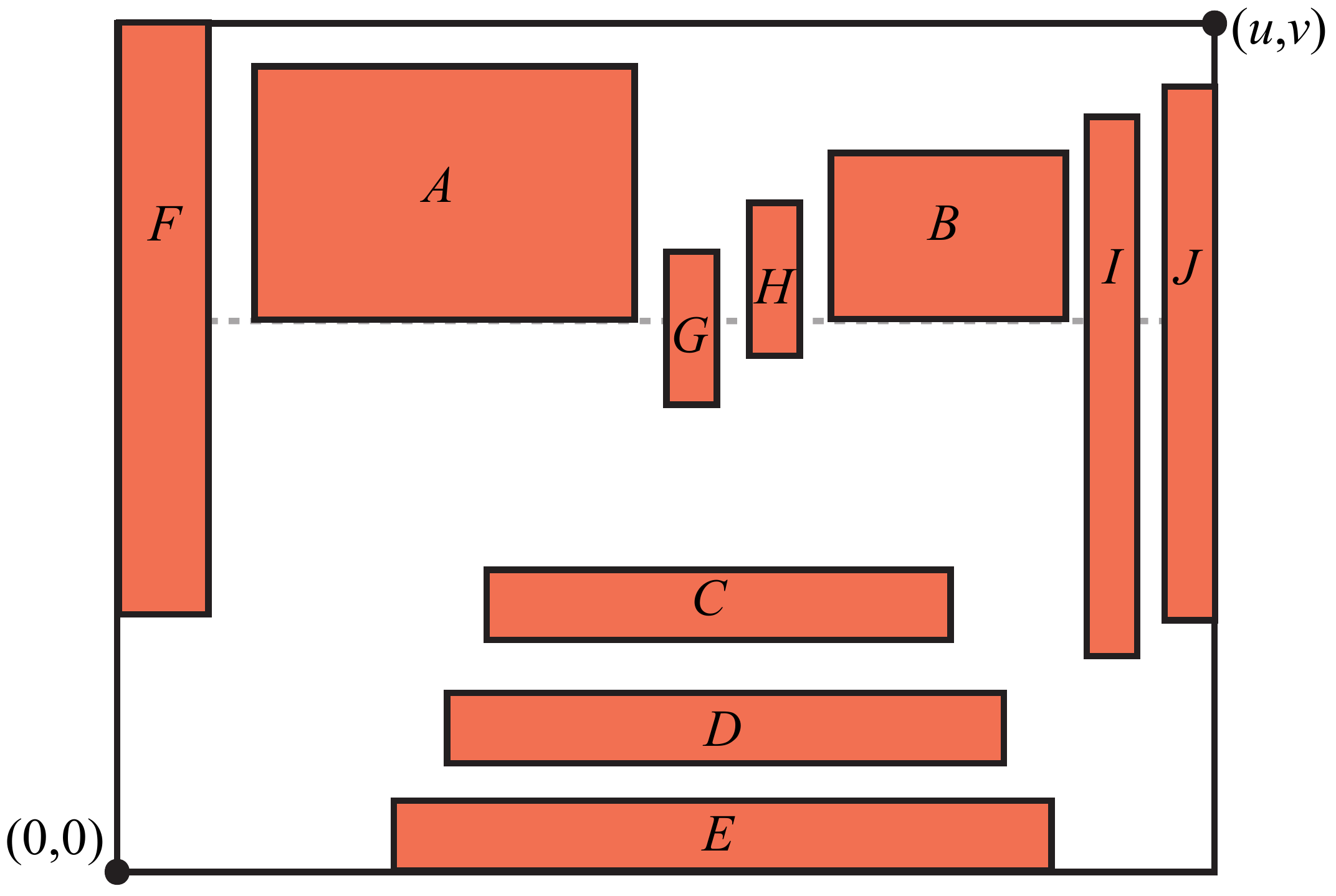}
    \caption{An example to illustrate the proof of Lemma~\ref{lonely}}
    \label{fig:lonely rectangle}
\end{figure}

\begin{lemma}\label{lonely}  Suppose $n\geq 6$ and let $S$ be any rectangle visibility representation of $K_n$. 
Then $|{\mathcal{T}}_S  | =1.$
\end{lemma}
 
\begin{proof} By way of contradiction, suppose $|{\mathcal{T}}_S  | \geq 2.$ Fix any distinct rectangles $A, B \in {\mathcal{T}}_S$, and note that ${\mathcal{N}}(A)$ and ${\mathcal{N}}(B)$ are empty, as illustrated by the example in Figure~\ref{fig:lonely rectangle}. Without loss of generality, assume that $x_1^A \leq x_1^B$ and $y_2^A \geq y_2^B$ (i.e., the taller rectangle is on the left). We observe the following:
\medskip
 
$\bullet \;$ 
${\mathcal{W}}(A)$ is empty: 
If $F \in {\mathcal{W}}(A)$, then $F$ cannot see $B.$
\medskip
 
$\bullet \;$ 
$|{\mathcal{E}}(B)| \leq 1$: 
Otherwise, fix $I$, $J \in {\mathcal{E}}(B)$ with $x_1^I < x_1^J$, and note that $y_1^I, y_1^J\le y_1^B$. Then $y_2^I > y_2^B$ since $I$ must see $A$. But now $J$ cannot see $B$, a contradiction.
\medskip 
 
$\bullet \;$ 
${\mathcal{S}}(A)={\mathcal{S}}(B)$: 
If $C \in {\mathcal{S}}(A)$ then $y_2^C \leq y_1^A = y_1^B$. But $C$ sees $B$, so $C \in {\mathcal{S}}(B)$. Thus ${\mathcal{S}}(A) \subseteq {\mathcal{S}}(B)$ and, similarly, ${\mathcal{S}}(B) \subseteq {\mathcal{S}}(A)$. 
\medskip
 
  
$\bullet \;$ 
$|{\mathcal{S}}(A)| >1$:  
Otherwise $|{\mathcal{E}}(A)| \geq 4$, since $n \geq 6$ and ${\mathcal{N}}(A)$ and ${\mathcal{W}}(A)$ are empty. 
Since $|{\mathcal{E}}(B)|\leq 1$, this implies $|{\mathcal{E}}(A) \cap {\mathcal{W}}(B)| \geq 2$. Fix distinct $G$, $H \in {\mathcal{E}}(A) \cap {\mathcal{W}}(B)$ with $x_1^G  \leq x_1^H$.  Now $y_2^H > y_2^G$ since $H$ sees $A$. But then $G$ cannot see $B$, a contradiction.
\medskip
 
$\bullet \;$ 
$|{\mathcal{S}}(A)| \leq 1$:  
Otherwise, fix distinct $C$,$D \in {\mathcal{S}}(A)$ with $y_1^C  \geq y_1^D$. 
  Since ${\mathcal{S}}(A)={\mathcal{S}}(B),$ both $C$, $D$ see $A$ and $B$ from below. Now ${\mathcal{E}}(A)\cap {\mathcal{W}}(B)$ is empty, since if $G \in {\mathcal{E}}(A) \cap {\mathcal{W}}(B)$, then $D$ cannot see $G$. Since ${\mathcal{N}}(A)$ and ${\mathcal{W}}(A)$ are empty and $n \geq 6$, it follows that
  $|{\mathcal{E}}(A)| \geq 3$ and thus ${\mathcal{E}}(B) \geq 2$, a contradiction.
\medskip

Having shown that  $|{\mathcal{S}}(A)| > 1$ and $|{\mathcal{S}}(A)| \leq 1$, we have arrived at a contradiction, and we conclude that $|{\mathcal{T}}_S  | =1.$
\end{proof}

\begin{figure}
    \includegraphics[height=1.9in]{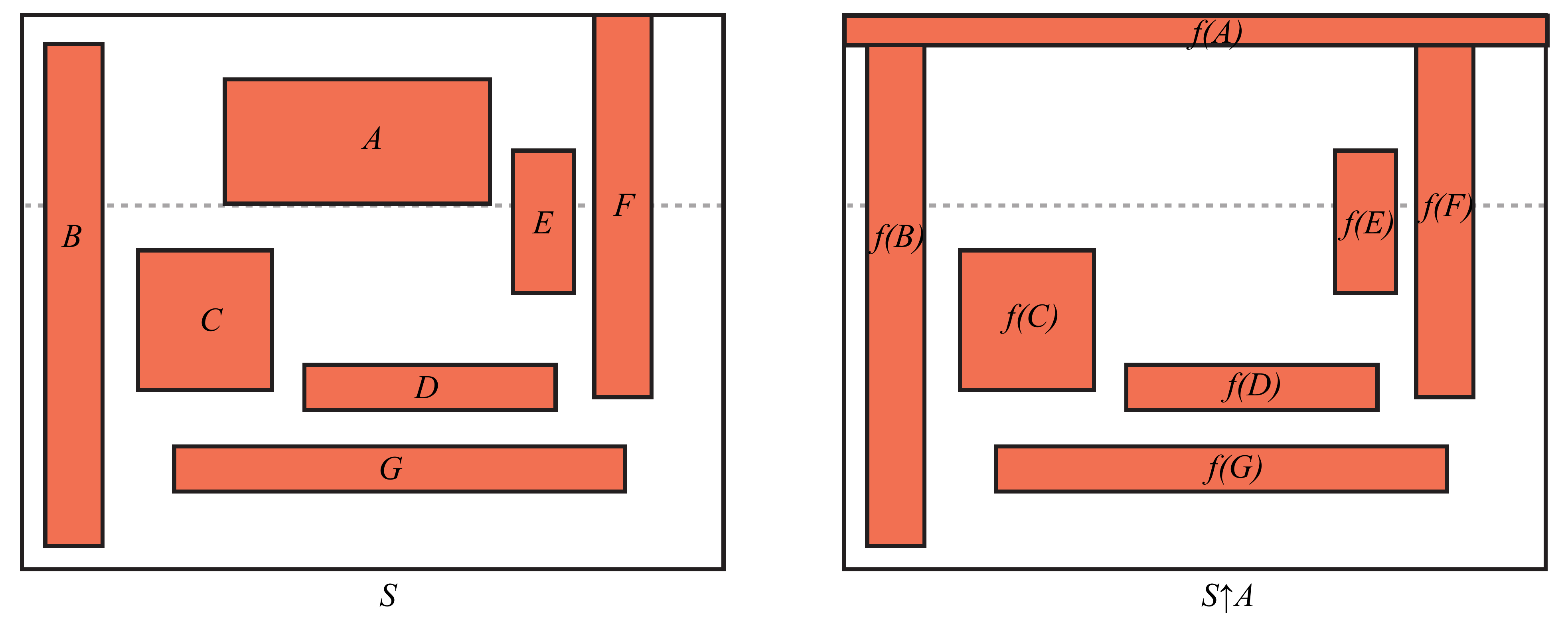}
    \caption{Applying the extracting operation $S \uparrow A$.}
    \label{fig:my_label}
\end{figure}

Another operation that will be useful is an extraction operation that can move a certain rectangle to the top row of the bounding box. Specifically, for a rectangle visibility representation $S$ of $K_n$ with bounding box $R = [0,u] \times [0,v]$ and a rectangle $A \in \mathcal{T}_S$, we define $S \uparrow A$ to be the set of rectangles in $R$ given by	$$ S \uparrow A = \{f(X) : X \in S\},$$
where $f(A) = [0,u] \times [v-1,v],$ and where, for every $X \not= A$, $$ f(X) = [x_1^X,x_2^X] \times [y_1^X, \min\{y_2^X,v-1\}].$$   As illustrated in Figure~\ref{fig:my_label}, the function $f$ maps the rectangle $A$ to the top row of $R$ and maps no other rectangle to that row. Furthermore, $f(X) \subseteq X$ for all $X \not= A$, so the rectangles of $ S \uparrow A$ do not overlap.

\begin{lemma}\label{extract} For any $n \geq 1$ let $S$ be a rectangle visibility representation of $K_n$ with bounding box $R = [0,u] \times [0,v]$. If ${\mathcal{T}}_S =\{A\}$ then $S \uparrow A$ also represents $K_n$ and has bounding box $R$.
\end{lemma}
 
 \begin{proof} Since $f$ bijectively maps the $n$ rectangles of $S$ to the $n$ rectangles of $S\uparrow A$, and since these representations share the same bounding box $R$, it remains only to show that $f$ preserves adjacency.  
 	
Since ${\mathcal{N}}(A)$ is empty and the graph is complete, $S$ is partitioned as $$ S = \{A\} \cup {\mathcal{E}}(A) \cup {\mathcal{W}}(A) \cup {\mathcal{S}}(A).$$ For any distinct $X$ and $Y$ in $S$, we claim  $f(X)$ and $f(Y)$ see each other: 

\smallskip

\textbf{Case 1. $X=A$ and $Y \in {\mathcal{E}}(A)$.} Since ${\mathcal{T}}_S=\{A\}$, note ${\mathcal{N}}(Y)$ is empty. So ${\mathcal{N}}(f(Y))=\{f(A)\}$, and $f(Y)$ sees $f(X)$ vertically. 

\smallskip

\textbf{Case 2. $X=A$ and $Y \in {\mathcal{S}}(A)$.} Note 
${\mathcal{S}}(f(A)) \supseteq {\mathcal{S}}(A)$ 
and so  $f(Y)$ sees $f(X)$ vertically.

\smallskip

\textbf{Case 3. $X\not=A$ and $Y \in {\mathcal{S}}(A)$.} Since ${\mathcal{T}}_S=\{A\}$, any line of sight between $X$ and $Y$ must be contained in the region below the top row of $R$. The only change to this region in $S \uparrow A$ is the removal of $A$, so $f(X)$ still sees $f(Y)$ along the original line of sight.

\smallskip

\textbf{Case 4. $\{X,Y\} \subseteq {\mathcal{E}}(A)$.} If $X$ sees $Y$ vertically, say with $X$ above $Y$, then $y_1(X) > y_1(A)$ so that $Y$ can see $A$.  Since ${\mathcal{T}}_S=\{A\}$, $X$ must see $Y$ horizontally. If the only line of visibility from $X$ to $Y$ were in the top row of $R$, then $y_2^X =y_2^Y = v$ and one of $X$,$Y$ could not see $A$. Therefore, $X$ must see $Y$ in a lower row, and so $f(X)$ still sees $f(Y)$ horizontally.

\smallskip

\textbf{Case 5. $X \in {\mathcal{E}}(A)$ and $Y \in {\mathcal{W}}(A)$.} Since ${\mathcal{T}}_S=\{A\}$, $X$ must see $Y$ horizontally. If $X$ sees $Y$ in any row below the top row of $R$, then $f(X)$ still sees $f(Y)$ in that same row. If $X$ {\em only} sees $Y$ in the top row of $R$, then $y_2^A < v$. Since ${\mathcal{T}}_S=\{A\}$, both $X$,$Y$ see $A$ horizontally in the top row of $A$. Since the bottom row of $f(A)$ is above the top row of $A$, now $f(X)$ sees $f(Y)$ horizontally in what was the top row of $A$.

\smallskip

By symmetry, these cover all possible cases.
\end{proof}

\begin{figure}[ht!]
 \centerline{\includegraphics[height=1.9in]{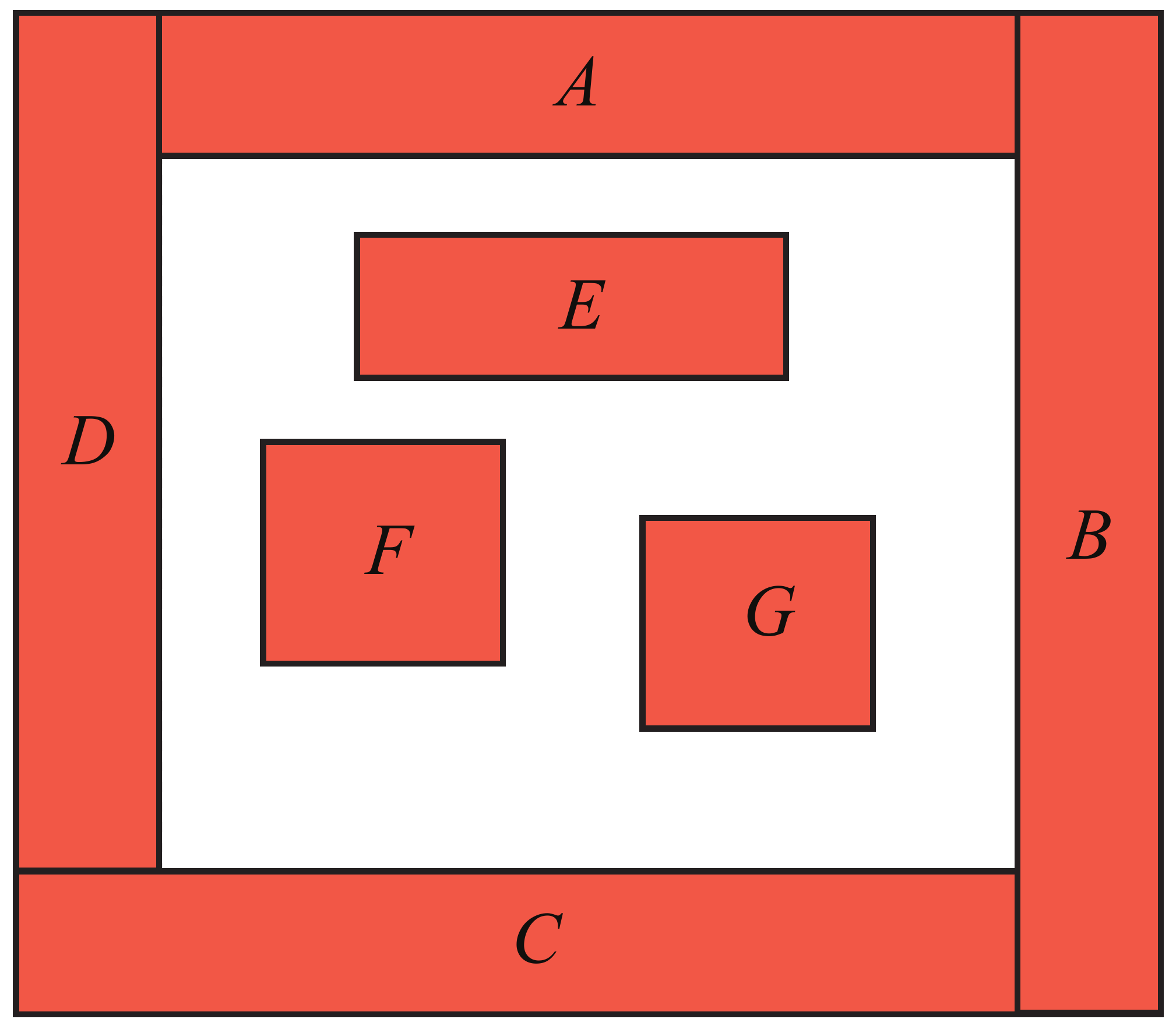} }
 \caption{The seven rectangles in the proof of Theorem~\ref{htk7}.}
\label{fig:K7proof}
 \end{figure}
 
\begin{lemma}\label{4out} Assume $n\geq 6$ and $K_n$ has a rectangle visibility representation with bounding box $R = [0,u] \times [0,v]$. Then $K_n$ has a rectangle visibility representation with bounding box $R$ in which the boundary of $R$ is covered by 4 rectangles of height or width 1.
\end{lemma}

\begin{proof}
	Apply 
	Lemma \ref{extract} successively in each of the four directions. Each time, a rectangle is brought to the corresponding boundary without changing the bounding box $R$.
\end{proof}

Recall that a \textit{\textbf{bar visibility graph}} $G$ is a graph representable with a set of disjoint horizontal bars in the plane, with edges between bars that have vertical lines of sight between them.  All bar visibility graphs are planar \cite{Duchet83, Tamassia86, Wismath85}.

\begin{lemma}
Suppose $S$ is an RV-representation for a graph $G$.  If we partition the edges of $G$ into those with vertical lines of sight and horizontal lines of sight in $S$, then the subgraphs $G_V(S)$ and $G_H(S)$ of $G$ with these edges are bar visibility graphs, and hence planar graphs. \label{bar-planar}
\end{lemma}

\begin{proof}
Replace each rectangle $A$ in $S$ with a horizontal line segment at the top edge of $A$.  This is a bar visibility representation of $G_V$.  Rotating $S$ by 90 degrees and then replacing each rectangle by its new top edge yields a bar visibility representation of $G_H$.
\end{proof}

\begin{figure}
\centering
\includegraphics[width=.3\textwidth]{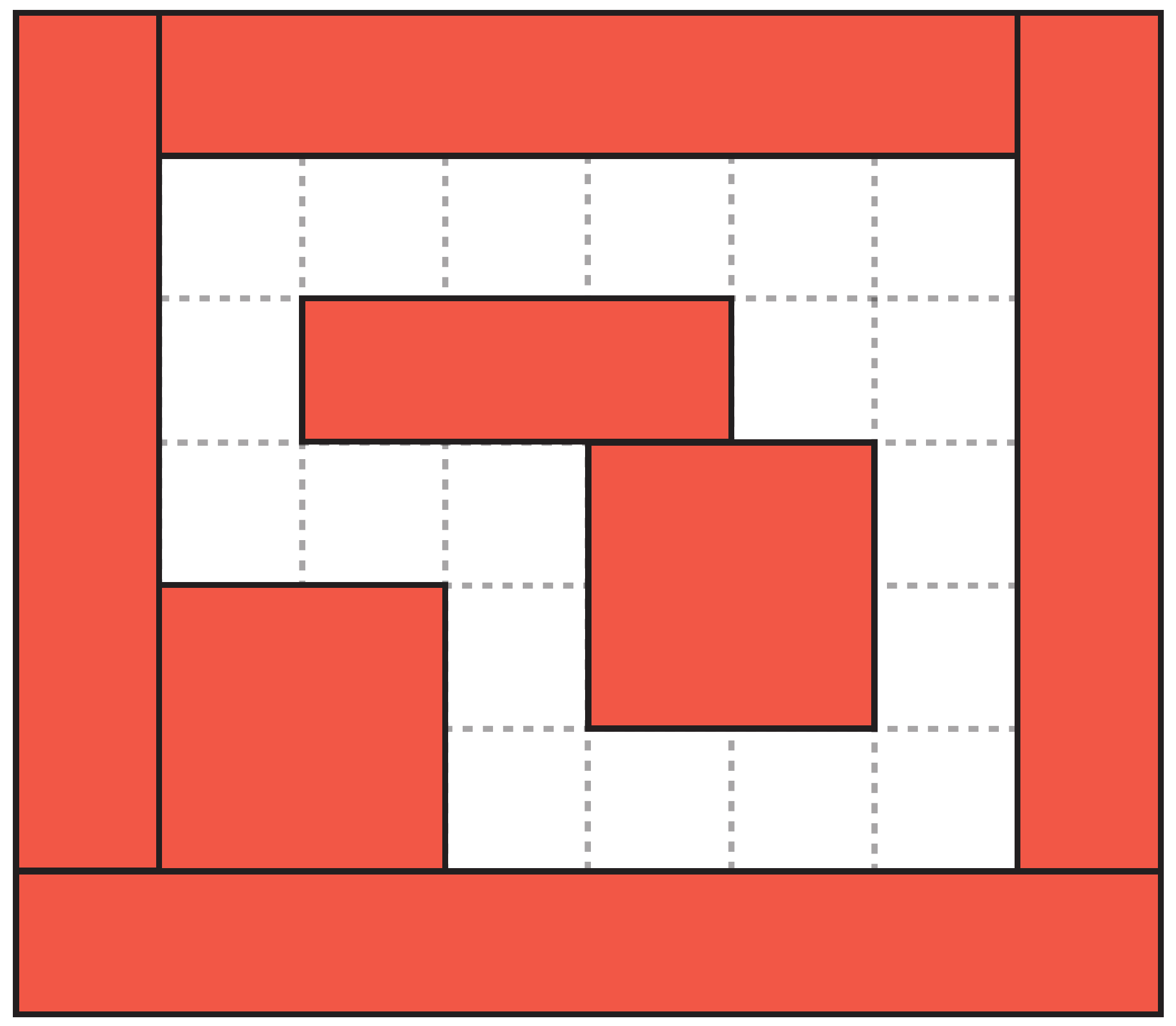} \hspace{1cm} \includegraphics[width=.4\textwidth]{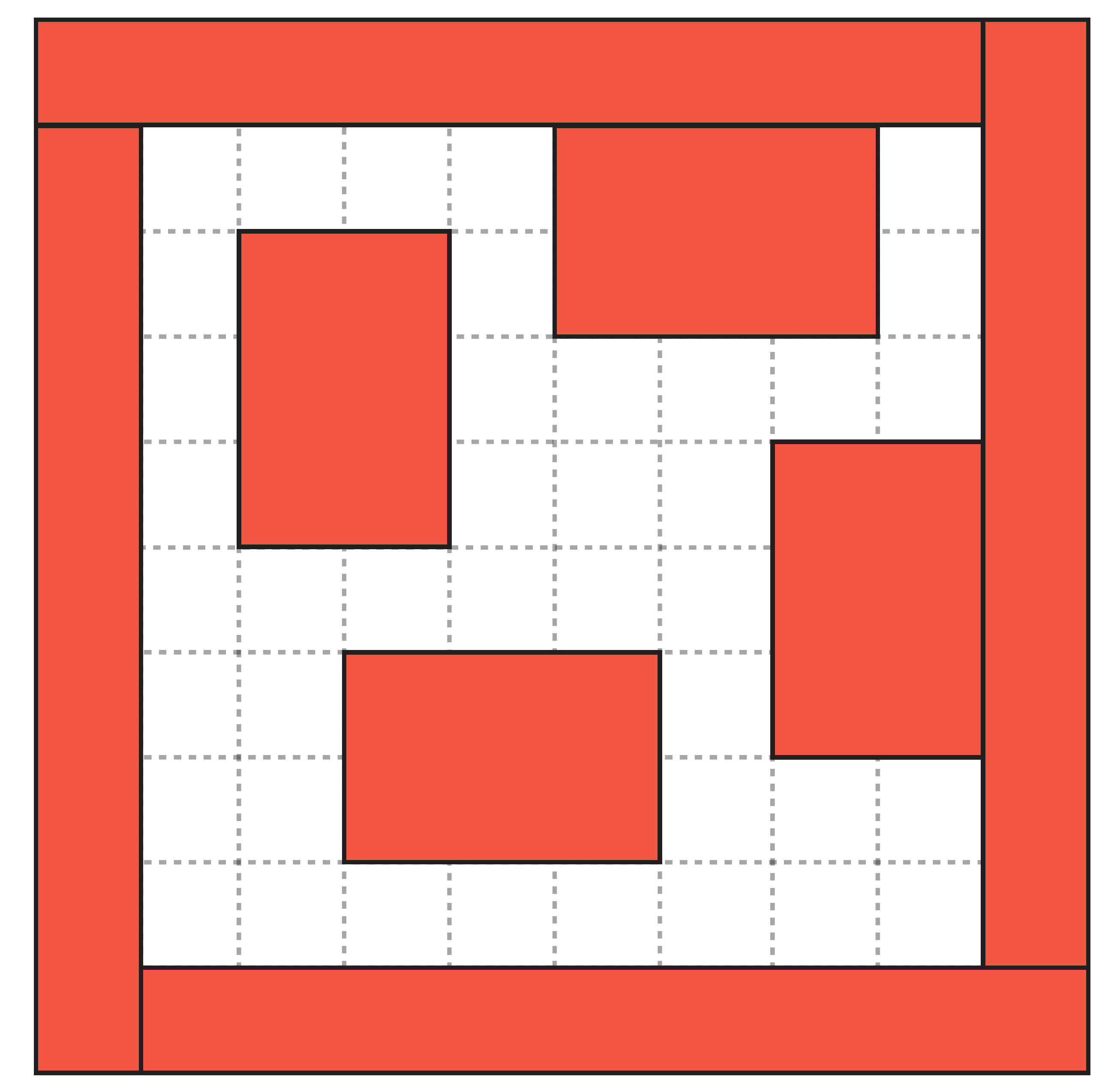}
\caption{A representation of $K_7$ in a $7 \times 8$ bounding box, and a representation of $K_8$ in a $10 \times 10$ bounding box.}
\label{fig:complete-graph-examples}
\end{figure}

\begin{theorem} \label{htk7} The complete graph $K_7$ has $\height(K_7)=7$.
\end{theorem}
\begin{proof}
Suppose $S$ is an RV-representation of $K_7$ with minimum height. By Lemma \ref{4out}, we may assume the boundary of the bounding box $R$ is covered by 4 rectangles of height or width 1, as in Figure \ref{fig:K7proof}.  Label these $A$, $B$, $C$, and $D$ clockwise from the top.

The remaining 3 rectangles $E$, $F$, and $G$ in the interior induce a 3-clique. If the 3 edges of this clique all correspond to horizontal lines of sight, then, together with rectangles $B$ and $D$, the edges of a 5-clique are represented entirely by horizontal lines of sight.  This is a 5-clique in $G_H$, which is impossible by Lemma~\ref{bar-planar}.  Similarly, the edges among $E$, $F$, and $G$ cannot all be vertical lines of sight. Rotating and renaming if necessary, assume $E$ sees $F$ and $G$ vertically and $F$ sees $G$ horizontally. Then $F$ and $G$ must be on the same side of $E$ and we may assume $F,G \in {\mathcal{S}}(E)$ and $G \in {\mathcal{E}}(F)$, as shown in Figure~\ref{fig:K7proof}. 

The following five horizontal lines of sight must occupy five distinct rows in $R$:
 $BD$, $BE$,  $FG$, $BF$, and $DG$. To see why, notice first that edge $BD$ must have its own row to reach all the way across the representation. The row for $BE$ must be distinct from $FG$, $BF$, and $DG$ since $F,G \in {\mathcal{S}}(E)$. Rows for $FG$ and $BF$ are distinct since $B,G \in {\mathcal{E}}(F)$. Rows for $FG$ and $DG$ are distinct since $F, D \in {\mathcal{W}}(G)$. Rows for $BF$ and $DG$ are distinct since $G \in {\mathcal{E}}(F)$.
 
Since $A$ and $C$ each take their own row by construction, it follows that $R$ has height at least 7. The representation of $K_7$ in Figure \ref{fig:complete-graph-examples} proves equality.
\end{proof}

\begin{theorem}\label{wk7} The complete graph $K_7$ has $\width(K_7)=8$.
\end{theorem}
\begin{proof}  If $\width(K_7)=7$, then Lemma~\ref{lowerbds}(iii) would guarantee an RV-representation $S$ in a $7 \times 7$ box $R$.
By Lemma \ref{4out} and Theorem \ref{htk7}, we may assume the boundary of $R$ is covered by 4 rectangles of height or width 1. Label these $A,B,C,D$ clockwise from the top.  As before, assume $F,G \in {\mathcal{S}}(E)$ and $G \in {\mathcal{E}}(F)$. 

The following six vertical lines of sight must occupy six distinct columns in $R$:
 $AC,AF,EF,EG,AG,CE$. To see why, note that edge $AC$ must have its own column to reach all the way across the representation. Any column meeting $F$ cannot meet $G$, since $G \in {\mathcal{E}}(F)$. Columns for $CE$,$EF$,$EG$ are distinct since $C,F,G \in {\mathcal{S}}(E)$. Columns for $EG$,$AG$ are distinct since $A,E \in {\mathcal{N}}(G)$. Columns for $AF$,$EF$ are distinct since $A,E \in {\mathcal{N}}(F)$. Columns for $AF$,$CE$ are distinct since $E \in {\mathcal{N}}(F)$. Columns for $AG$,$CE$ are distinct since $E \in {\mathcal{N}}(G)$.
 
 Thus $R$ has width at least 8.  By Figure \ref{fig:complete-graph-examples}, $\text{width}(K_7)=8$. 
\end{proof}

\begin{corollary}
$\area(K_7)=56$ and $\perimeter(K_7)=30$.
\end{corollary}

\begin{proof}
Theorems~\ref{htk7} and \ref{wk7} show that any representation of $K_7$ requires height at least 7 and width at least 8.  Figure~\ref{fig:complete-graph-examples} shows a representation of $K_7$ with height 7 and width 8 exactly.  Therefore this representation also has smallest area and perimeter by Lemma~\ref{lemma-same-height-width}.
\end{proof}

\begin{theorem} \label{htk8}
The complete graph $K_8$ has $\rm{height}(K_8)=10$.
\end{theorem} 

\begin{proof} Suppose $S$ is an RV-representation of $K_8$ with minimum height. By Lemma \ref{4out}, we may assume the boundary of the bounding box $R$ is covered by 4 rectangles of height or width 1, as in Figure \ref{4out}. Label these $A$, $B$, $C$, and $D$ clockwise from the top.

Rectangles $E$, $F$, $G$, and $H$ in the interior must induce a 4-clique. By Lemma~\ref{bar-planar}, the edges in this clique that correspond to vertical lines of sight must form a triangle-free subgraph (and similarly for the horizontal edges).  Up to isomorphism, there are only two decompositions of $K_4$ into a pair of triangle-free graphs, shown in Figure~\ref{fig:K8proof}.
 
 \begin{figure}[ht!]
 \centerline{\includegraphics[height=1in]{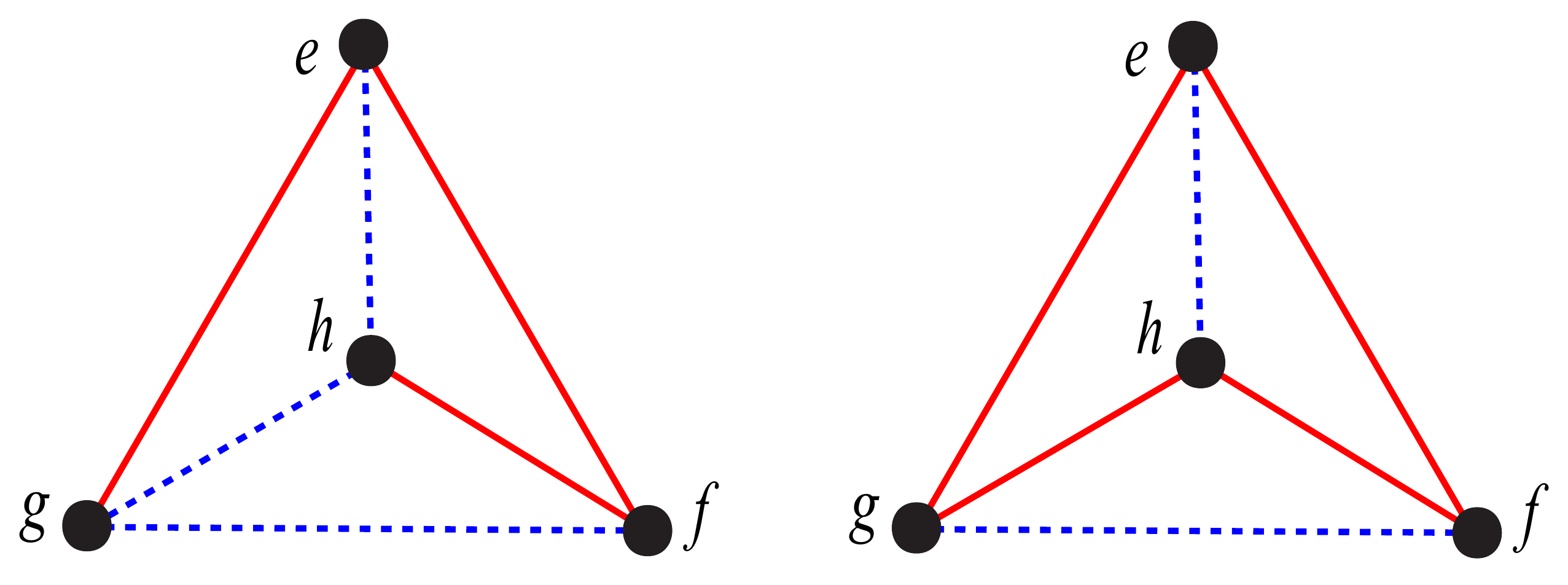}}
 \caption{The partition of the edges of $K_8$ in Theorem~\ref{htk8}.}
\label{fig:K8proof} 
 \end{figure}

First we claim that the graph on the right of Figure~\ref{fig:K8proof} is impossible. To see why, let the solid edges denote vertical lines of sight. Since $FG$ is horizontal and $EF$ and $EG$ are vertical, $F$ and $G$ are on the same side of $E$. Assume $F,G \in {\mathcal{S}}(E)$. Since $EH$ is horizontal and $FH$ and $GH$ are vertical, $F, G \in {\mathcal{S}}(H)$. Renaming if necessary, $E \in {\mathcal{E}}(H)$ and $F \in {\mathcal{E}}(G)$. So $x_2^G \leq x_1^F$. But $x_1^E < x_2^G$ and $x_1^F < x_2^H$. It follows that $x_1^E < x_2^H$, contradicting that $E \in {\mathcal{E}}(H)$. A similar argument elimates the case when the solid edges denote horizontal lines of sight.

Next we claim that the graph on the left requires $R$ to contain at least 10 rows. The graph is symmetric in solid and dotted edges, so assume the solid edges denote vertical lines of sight. Since $FG$ is horizontal and $EF$ and $EG$ are vertical, $F$ and $G$ are on the same side of $E$. Assume $F,G \in {\mathcal{S}}(E)$.
Since $EH$ is horizontal and $EF$ and $FH$ are vertical, $E$ and $H$ are on the same side of $F$. So $H \in {\mathcal{N}}(F)$.
Since $FH$ is vertical and $GH$ and $FG$ are horizontal, $F$ and $H$ are on the same side of $G$. Assume $F,H \in {\mathcal{W}}(G)$.
Since $EG$ is vertical and $GH$ and $EH$ are horizontal, $E$ and $G$ are on the same side of $H$. So $H \in {\mathcal{W}}(E)$. 

The following 8 horizontal lines of sight occupy distinct rows in $R$: 
 $$BD, DE, EH, BH, GH, DG, FG, BF. $$
Edge $BD$ must have its own row to reach all the way across the representation. Any row meeting $E$ cannot meet $F$ or $G$, since $F,G \in {\mathcal{S}}(E)$. Any row meeting $F$ cannot meet $H$, since $H \in {\mathcal{N}}(F)$. 
Rows $DE$ and $EH$ are distinct since $D,H \in {\mathcal{W}}(E)$.
Rows $EH$ and $BH$ are distinct since $B,E \in {\mathcal{E}}(H)$.
Rows $BH$ and $GH$ are distinct since $B,G \in {\mathcal{E}}(H)$.
Rows $GH$ and $DG$ are distinct since $D,H \in {\mathcal{W}}(G)$.
Rows $DG$ and $FG$ are distinct since $F,D \in {\mathcal{W}}(G)$.
Rows $FG$ and $BF$ are distinct since $B,G \in {\mathcal{E}}(F)$.
Rows $DE$ and $BH$ are distinct since $H \in {\mathcal{W}}(E)$.
Rows $BH$ and $DG$ are distinct since $H \in {\mathcal{W}}(G)$.
Rows $DG$ and $BF$ are distinct since $F \in {\mathcal{W}}(G)$.
 
Since $A$ and $C$ each take their own row by construction, $R$ has height $\geq 10$. The representation of $K_8$ shown in Figure~\ref{fig:complete-graph-examples} proves equality.
\end{proof}

\begin{theorem} \label{wk8} 
The complete graph $K_8$ has $\width(K_8)=10$.
\end{theorem} 

\begin{proof} Since, by definition, $\width(K_8) \geq \rm{height}(K_8) =10$, the representation of $K_8$ shown in Figure \ref{fig:complete-graph-examples} proves equality.
\end{proof}

\begin{corollary}
The complete graph $K_8$ has $\area(K_8)=100$ and $\perimeter(K_8)=40$.
\end{corollary}

\begin{proof}
Theorems~\ref{htk8} and \ref{wk8} show that any representation of $K_8$ requires height at least 10 and width at least 10.  Figure~\ref{fig:complete-graph-examples} shows a representation of $K_8$ with height 10 and width 10 exactly.  Therefore this representation also has smallest area and perimeter by Lemma~\ref{lemma-same-height-width}.
\end{proof}

Note that $K_8$ is the largest complete RVG \cite{hutchinson1999representations}, so we can't investigate the size of RV-representations of larger complete graphs.  But using the disjoint unions of complete graphs and Lemma~\ref{NewDisjointLem}, we can construct graphs on $n$ vertices whose RV-representations are larger than the empty graph for all $n \geq 8$, as follows.

\begin{corollary} \label{disjoint complete graphs}
Fix any positive integer $n$ and write $n=8q+r$ for integers $q$ and $r$ with  $0 \leq r < 8$.  Define the graph $G_n = q K_8 + E_r$, which has $q$ disjoint copies of $K_8$ and $r$ isolated vertices.  Then $G_n$ has $n$ vertices, $\height(G_n) = \width(G_n)=n+2q$, $\area(G_n)=(n+2q)^2$ and $\perimeter(G_n)=4n+8q$.
\end{corollary}

\begin{proof}
By Lemma~\ref{NewDisjointLem}, $\height(G_n)= q \cdot \height(K_8)+ \height(E_r)$.  By Theorem~\ref{htk8}, $\height(K_8)=10$, and again by Lemma~\ref{NewDisjointLem}, $\height(E_r)=r$.  So $\height(G_n)=10q+r=n+2q$.

Since the bounding boxes of the representations of $K_8$ and $E_r$ are square, the same argument holds for $\width(G_n)$.  Since the representations of $G_n$ with minimum height and width are the same, these representations also yield the minimum area and perimeter of $G_n$ by Lemma~\ref{lemma-same-height-width}.
\end{proof}

\begin{corollary}
Among all rectangle visibility graphs with $n \geq 7$ vertices, the empty graph $E_n$ does not have the largest width, area, or perimeter.
Among all rectangle visibility graphs with $n \geq 8$ vertices, the empty graph $E_n$ does not have the largest height.
\end{corollary}


\section{Directions for Further Research} \label{Conclusions}

We conclude with a number of open problems and questions that could further this line of research.

\begin{enumerate}[wide]
\item We have established that for $n=7,8$ the complete graph exceeds the empty graph in area, perimeter, height, and width. We have also shown that for $n>8$, the empty graph does not maximize any of these parameters. Accordingly, it is natural to ask, in general, which rectangle visibility graph(s) with $n$ vertices have largest height, perimeter, width, and area?  Note that for $n>8$, $K_n$ is not a rectangle visibility graph \cite{hutchinson1999representations}.  Furthermore, when $r=7$ in Corollary~\ref{disjoint complete graphs}, we can replace $E_7$ by $K_7$ to obtain a graph with width $n+2q+1$. Are there any other graphs that beat the graph in Corollary~\ref{disjoint complete graphs}?

\item Tables~\ref{table-separating1} and~\ref{table-separating-new} show that, for each pair of parameters of area, perimeter, height, and width, there are pairs of graphs that share the same number of vertices, but that are ordered oppositely by that pair of parameters.  However, it does not consider triples and quadruples of parameters.  For example, is there a pair of graphs $G_1$ and $G_2$ for which $\area(G_1)<\area(G_2)$ but both $\height(G_2)<\height(G_1)$ and $\width(G_2)<\width(G_1)$?

\item Tables~\ref{table-separating2} and~\ref{table-separating2b} show that, for each pair of parameters of area, perimeter, height, and width, there is a graph with two representations that are ordered oppositely by that pair of parameters.  However, it does not consider triples and quadruples of parameters.  For example, is there a graph $G$ that requires three distinct RV-representations to minimize its area, perimeter, and height?

\item Say that an RV-representation $S$ is \textit{compressible} if we can delete a row or column of $S$ and still have a representation of the same graph. For a given number of vertices $n$, which graphs have the largest incompressible representations, in terms of area, perimeter, height, or width?  How large are these values?

\item We might consider additional measures of size in terms of the rectangles in an RV-representation, rather than the bounding box.  For example, for an RV-representation $S$, say $\recarea(S)$ is the area of the largest rectangle in $S$.  Then $\recarea(G)$ is the smallest value of $\recarea(S)$ for any RV-representation of $G$.  Which graphs $G$ with $n$ vertices have largest $\recarea(G)$?

\item We conjecture that if $\height(G)=2$ then $G$ is outerplanar.  Is this true?  Can we characterize other families of graphs of specific area, perimeter, height, or width?

\item We can consider other dimensions. For dimension 1, what is the minimum length of an integer bar visibility graph on $n$ vertices?  For 3-dimensional box visibility graphs, there are many parameters measuring the size of an integer box visibility representation.  Which graphs require the largest 3-dimensional representation, as measured by these parameters?  Note that in \cite{Fekete99} Fekete and Meijer proved that $K_{56}$ is a 3-dimensional box visibility graph.
 
 \end{enumerate}

\bibliographystyle{plain}
\bibliography{RVGbib}
\end{document}